\documentclass{article}
\usepackage{color}
\usepackage{amsmath,amssymb}
\usepackage[titletoc]{appendix}
\catcode`\@=11 \@addtoreset{equation}{section}

\catcode`\@=12

\newtheorem{Theorem}{Theorem}[section]
\newtheorem{Proposition}{Proposition}[section]
\newtheorem{Lemma}{Lemma}[section]
\newtheorem{Remark}{Remark}[section]
\newtheorem{Corollary}{Corollary}[section]
\newtheorem{Definition}{Definition}[section]

\newcommand{\bTheorem}[1]{
\begin{Theorem} \label{T#1} }
\newcommand{\eT}{\end{Theorem}}

\newcommand{\bProposition}[1]{
\begin{Proposition} \label{P#1}}
\newcommand{\eP}{\end{Proposition}}

\newcommand{\bLemma}[1]{
\begin{Lemma} \label{L#1} }
\newcommand{\eL}{\end{Lemma}}

\newcommand{\bCorollary}[1]{
\begin{Corollary} \label{C#1} }
\newcommand{\eC}{\end{Corollary}}

\newcommand{\bFormula}[1]{
\begin{equation} \label{#1}}
\newcommand{\eF}{\end{equation}}

\newcommand{\lpq}{L_p(0,T;L_q(\Omega_0))}

\newcommand{\vr}{\varrho}

\newcommand{\vrt}{\tilde{\vr}}

\newcommand{\vut}{\tilde{\vu}}

\newcommand{\expnorm}{\|e^{\gamma t}(\sigma,\vv)\|_{\dot{\cal Y}(\infty)}}
\newcommand{\expnormzw}{\|e^{\gamma t}(z,\vw)\|_{\dot{\cal Y}(\infty)}}
\newcommand{\bX}{\vc{X}}
\newcommand{\bY}{\vc{Y}}

\newcommand{\vu}{\vc{u}}

\newcommand{\vv}{\vc{v}}
\newcommand{\bv}{\vc{v}}

\newcommand{\vw}{\vc{w}}

\newcommand{\vf}{\vc{f}}
\newcommand{\vF}{\vc{F}}
\newcommand{\bF}{\vc{F}}

\newcommand{\vV}{\vc{V}}
\newcommand{\vVt}{\tilde{\vc{V}}}

\newcommand{\vR}{\vc{R}}
\newcommand{\bR}{\vc{R}}
\newcommand{\vE}{\vc{E}}
\newcommand{\vc}[1]{{\bf #1}}
\newcommand{\qed}{\rightline{ $\square$}}
\renewcommand{\div}{{\rm div}}
\newcommand{\dv}{{\rm div}}
\newcommand{\Div}{{\rm div}_x}
\newcommand{\Grad}{\nabla_x}

\newcommand{\tn}[1]{\mbox {\F #1}}
\newcommand{\dx}{{\rm d} {x}}

\renewcommand{\d}{{\rm d} }
\newcommand{\dt}{{\rm d} t }
\newcommand{\ds}{{\rm d} s }

\newcommand{\bProof}{{\bf Proof: }}

\newcommand{\ep}{\varepsilon}
\newcommand{\de}{\partial}

\renewcommand{\div}{{\rm div}\,}

\newcommand{\bXu}{\bX_{\vu}}
\newcommand{\etab}{\bar \eta}
\newcommand{\vvb}{\bar \vv}

\newcommand{\CL}{{\mathcal L}}
\newcommand{\CC}{{\mathcal C}}
\newcommand{\CR}{{\mathcal R}}
\newcommand{\CF}{{\mathcal F}}
\newcommand{\CD}{{\mathcal D}}
\newcommand{\CS}{{\mathcal S}}
\newcommand{\CA}{{\mathcal A}}
\newcommand{\CB}{{\mathcal B}}
\newcommand{\CT}{{\mathcal T}}
\newcommand{\CH}{{\mathcal H}}
\newcommand{\BN}{{\mathbb T}}

\newcommand{\ku}{{\bf k}_{\vut}}
\newcommand{\kv}{{\bf k}_{\vv}}
\newcommand{\kw}{{\bf k}_{\vw}}
\newcommand{\pd}{\partial}
\newcommand{\bI}{\vc{I}}
\newcommand{\lr}[1]{\left( #1 \right)}
\newcommand{\I}{{\bf I}}
\newcommand{\eq}[1]{\begin{equation}
\begin{split}
#1
\end{split}
\end{equation}}

\font\F=msbm10 scaled 1000
\newcommand{\R}{\mathbb{R}}
\newcommand{\C}{\mbox{\F C}}

\newcommand{\N}{\mathbb{N}}

\definecolor{grey}{rgb}{0.85,0.85,0.85}
\date{}

\long\def\greybox#1{%
    \newbox\contentbox%
    \newbox\bkgdbox%
    \setbox\contentbox\hbox to \hsize{%
        \vtop{
            \kern\columnsep
            \hbox to \hsize{%
                \kern\columnsep%
                \advance\hsize by -2\columnsep%
                \setlength{\textwidth}{\hsize}%
                \vbox{
                    \parskip=\baselineskip
                    \parindent=0bp
                    #1
                }%
                \kern\columnsep%
            }%
            \kern\columnsep%
        }%
    }%
    \setbox\bkgdbox\vbox{
        \color{grey}
        \hrule width  \wd\contentbox %
               height \ht\contentbox %
               depth  \dp\contentbox
        \color{black}
    }%
    \wd\bkgdbox=0bp%
    \vbox{\hbox to \hsize{\box\bkgdbox\box\contentbox}}%
    \vskip\baselineskip%
}

\begin{document}


\title{Compressible Navier-Stokes system on a moving domain in the $L_p-L_q$ framework}

\author{
Ond{\v r}ej Kreml $^1$
\and {\v S}{\' a}rka Ne{\v c}asov{\' a} $^1$
\and Tomasz Piasecki $^2$
}

\maketitle

\bigskip

\centerline{$^1$Institute of Mathematics of the Academy of Sciences of
the Czech Republic} \centerline{\v Zitn\' a 25, 115 67 Praha
Czech Republic}
\bigskip

\centerline{$^2$  Institute of Applied Mathematics and Mechanics, University of Warsaw}
\centerline{Banacha 2, 02-097 Warszawa, Poland}

\medskip

\begin{abstract}
We prove the local well-posedness for the barotropic compressible Navier-Stokes system on a moving domain, a motion of which is determined by a given vector field $\vV$, in a maximal $L_p-L_q$ regularity framework. Under additional smallness assumptions on the data we show that our solution exists globally in time and satisfies a decay estimate. In particular, for the global well-posedness we do not require exponential decay or smallness of $\vV$ in $L_p(L_q)$. However, we require exponential decay and smallness of its derivatives. 

\end{abstract}

\medskip


\medskip


\section{Introduction}\label{i}

We consider a barotropic flow of a compressible viscous fluid in the absence of external forces described by the isentropic compressible Navier-Stokes system
\bFormula{i1a}
\partial_t \vr + \Div (\vr \vu) = 0,
\eF
\bFormula{i1b}
\partial_t (\vr \vu) + \Div (\vr \vu \otimes \vu) + \Grad \pi(\vr) = \Div \tn{S}(\Grad \vu),
\eF
where $\vr$ is the density of the fluid and $\vu$ denotes the velocity. 
We assume that the stress tensor $\tn{S}$ is determined by the standard Newton rheological law
\bFormula{i4}
\tn{S} (\Grad \vu) = \mu \left( \Grad \vu +
\Grad^t \vu - \frac{2}{3} \Div \vu \tn{I} \right) + \zeta \Div \vu \tn{I}
\eF
with constant viscosity coefficients $\mu > 0$ and $\zeta \geq 0$. 
The pressure $\pi(\vr)$ is a given sufficiently smooth, nondecreasing function of the density. 
We assume the fluid occupies a time-dependent bounded domain $\Omega_t$, the motion of which is described by means of a given velocity field $\vc{V}(t,x)$, 
where $t \geq 0$ and
$x \in \mathbb{R}^3$. More precisely, we assume that if $\vc{X}$ solves the following system of ordinary differential equations
\begin{equation*}
\frac{{\rm d}}{{\rm d}t} \vc{X}(t, x) = \vc{V} \Big( t, \vc{X}(t, x) \Big),\ t > 0,\ \vc{X}(0, x) = x,
\end{equation*}
we set
\[
\Omega_\tau = \vc{X} \left( \tau, \Omega_0 \right),\]
where $\Omega_0 \subset \mathbb{R}^3$ is a given bounded domain at initial time $t = 0$. Moreover we denote $\Gamma_\tau = \partial \Omega_\tau$ and
\[ Q_\tau = \bigcup_{t \in (0,\tau)} \{t\}\times \Omega_t =: (0,\tau)\times \Omega_t.
\]
%
%
We consider the system \eqref{i1a}-\eqref{i1b} supplied with 
the Dirichlet boundary conditions 
\bFormula{i6a}
(\vu - \vc{V})|_{\Gamma_\tau} = 0 \ \mbox{for any}\ \tau \geq 0
\eF
%
%
and the initial conditions
\bFormula{i7}
\vr(0, \cdot) = \vr_0 ,\quad
\vu(0, \cdot) = \vu_0 \quad \mbox{in}\ \Omega_0.
\eF
%
%

The existence theory for system \eqref{i1a}-\eqref{i1b} on fixed domains is nowadays quite well developed. 
The existence of global weak solutions has been first established by Lions \cite{LI4}. 
This result has been later extended by Feireisl and coauthors (\cite{FNP}, \cite{EF61}, \cite{EF70}, \cite{EF71}) 
to cover larger class of pressure laws.
Strong solutions on fixed domains are known to exist locally in time or globally provided certain smallness assumptions 
on the data. For no-slip boundary conditions see among others \cite{MN1}, \cite{MN2}, \cite{V1}, \cite{V2} 
for the results in Hilbert spaces,
\cite{MZ1}, \cite{MZ2}, \cite{MZ3} in $L_p$ setting and \cite{ES} for a maximal $L_p-L_q$ regularity approach.   
The problem with slip boundary conditions on a fixed domain has been 
investigated by Zajaczkowski \cite{Za}, Hoff \cite{Ho} and, more recently, by Shibata and Murata \cite{Mur}, \cite{SM}  
in the $L_p - L_q$ maximal regularity setting. In \cite{PSZ},\cite{PSZ2} the approach from \cite{ES} has been adapted to treat a generalization of the compressible Navier-Stokes system describing a flow of a compressible mixture with cross-diffusion. For results on free boundary problems of the system \eqref{i1a}-\eqref{i1b},
we refer to \cite{Za1}, \cite{Za2} where the global existence of strong solutions in $L_2$ setting has been shown under the assumption that the domain is close to a ball and to \cite{S1} where a free boundary problem is treated in $L_p-L_q$ approach. 

The existence theory for the system \eqref{i1a}-\eqref{i1b} on a moving domain with a given motion of the boundary 
started to develop with the results for weak solutions obtained using a penalization method in \cite{FeNeSt} for 
no-slip bondary conditions and \cite{FKNNS} for slip conditions. These results have been recently generalized 
to the complete system with heat conductivity in \cite{KMNW1} and \cite{KMNW2}. 
The first weak-strong uniqueness 
result on a moving domain has been shown in \cite{dobo} in case of no-slip boundary condition. 
A generalization of this to slip conditions as well as a local existence result for strong solution 
for both types boundary conditions can be found in \cite{KNP}. There, the authors 
use the energy approach in $L_2$ setting for the existence result. 

The aim of this paper is to extend the existence theory for strong solutions on a moving domain 
to $L_p - L_q$ maximal regularity setting. We present a more detailed outline of the proof after stating our main result, however first let us resume the notation  used in the paper.

\subsection{Notation}
We use standard notations for Lebesgue spaces $L_p(\Omega_0)$ and Sobolev spaces $W^k_p(\Omega_0)$ with $k \in \N$ on a fixed domain $\Omega_0$. By $C_B(\Omega_0)$ we denote a space of bounded continuous functions on $\Omega_0$. 
Furthermore, for a Banach space $X$, $L_p(0,T;X)$ is a Bochner space of functions for which the norm
\begin{equation*}
\|f\|_{L_p(0,T;X)}=\left\{ 
\begin{array}{lr} 
\left( \int_0^T \|f(t)\|_X\, \dt\right)^{1/p} , & 1 \leq p < \infty,\\[3pt]
\textrm{ess sup}_{0 \leq t \leq T} \|f(t)\|_X, & p=\infty
\end{array}
\right.
\end{equation*}
is finite. Then 
$$
W^1_p(0,T;X)=\{f \in L_p(0,T;X): \; \de_t f \in L_p(0,T;X) \}. 
$$
For $p \geq 1$ we denote by $p'$ its dual exponent, i.e. $\frac 1p + \frac{1}{p'} = 1$. Next, we recall that for $0<s<\infty$ and $m$ a smallest integer larger than $s$ we define Besov spaces on domains as intermediate spaces
\begin{equation} \label{def:bsqp0} 
B^{s}_{q,p}(\Omega_0)=(L_q(\Omega_0),W^m_q(\Omega_0))_{s/m,p},
\end{equation}
where $(\cdot,\cdot)_{s/m,p}$ is the real interpolation functor, see \cite[Chapter 7]{Ad}. In particular,
\begin{equation} \label{def:bsqp}
B^{2(1-1/p)}_{q,p}(\Omega_0)=(L_q(\Omega_0),W^2_q(\Omega_0))_{1-1/p,p}=(W^2_q(\Omega_0),L_q(\Omega_0))_{1/p,p}.    \end{equation}
We shall not distinguish between notation of spaces for scalar and vector valued functions, i.e. we write $L_q(\Omega_0)$ instead of $L_q(\Omega_0)^3$  etc. However, we write vector valued functions in boldface.

For $p, q\in[1,\infty]$ we introduce the vector space
$L_p(0,T,L_q(\Omega_t))$ which consists of functions $f:Q_T\mapsto \mathbb R$ such that 
$$
t\mapsto \|f(t,\cdot)\|_{L_q(\Omega_t)}
$$
is measurable and $L_p$ integrable on time interval $(0,T)$, i.e. if
\begin{equation}
    \|f\|_{L_p(0,T,L_q(\Omega_t))} := \left\| \|f(t,\cdot)\|_{L_q(\Omega_t)} \right\|_{L_p(0,T)} \leq C.
\end{equation}
Similarly, for $p, q\in[1,\infty]$ and $k \in \N$  we introduce a vector space $L_p(0,T,W^k_q(\Omega_t))$ which consists of functions $f:Q_T\mapsto \mathbb R$ such that 
\begin{equation}
    \|f\|_{L_p(0,T,W^k_q(\Omega_t))} := \left\| \|f(t,\cdot)\|_{W^k_q(\Omega_t)} \right\|_{L_p(0,T)} \leq C.
\end{equation}


Let us also introduce a brief notation for the regularity class of the solution. Namely, for a function $g$ and a vector field $\vf$ defined on $(0,T)\times\Omega_t$ we define
\begin{equation} \label{def:X}
\|g,\vf\|_{{\cal{X}}(T)}=\|\vf\|_{L_p(0,T;W^2_q(\Omega_t))}+ \|\vf_t\|_{L_p(0,T;L_q(\Omega_t))}+\|g\|_{L_p(0,T;W^1_q(\Omega_t))}+\|g_t\|_{L_p(0,T;L_q(\Omega_t))}. 
\end{equation}
Next, for a pair $\tilde g,\tilde \vf$ defined of $(0,T) \times \Omega_0$ we define a norm  
\begin{equation} \label{def:Y}
\|\tilde g,\tilde \vf\|_{{\cal{Y}}(T)}=\|\tilde \vf\|_{L_p(0,T;W^2_q(\Omega_0))}+ \|\tilde \vf_t\|_{L_p(0,T;L_q(\Omega_0))}+\|\tilde g\|_{W^1_p(0,T;W^1_q(\Omega_0))}
\end{equation}
and a seminorm
\begin{equation} \label{def:dotY}
\|\tilde g,\tilde \vf\|_{\dot{\cal{Y}}(T)}=\|\tilde \vf\|_{L_p(0,T;W^2_q(\Omega_0))}+ \|\tilde \vf_t\|_{L_p(0,T;L_q(\Omega_0))}+\|\nabla \tilde g\|_{L_p(0,T;L_q(\Omega_0))}+\|\de_t \tilde g\|_{L_p(0,T;W^1_q(\Omega_0))}.
\end{equation}
Obviously, we denote by ${\cal X}(T)$ and ${\cal Y}(T)$ spaces for which the norms \eqref{def:X} and \eqref{def:Y}, respectively, are finite.
\begin{Remark}
Notice that the norm \eqref{def:Y} involves also $\|\tilde g_t\|_{L_p(W^1_q(\Omega_0))}$ while in \eqref{def:X} we have only $\|g_t\|_{L_p(L_q(\Omega_t))}$. The reason is that in the Lagrangian coordinates we are able to show higher regularity of the density, which does not correspond to equivalent regularity in Eulerian coordinates, see Section \ref{sec:equiv}.
\end{Remark}
Furthermore we denote 
\begin{align}
&\CH_q = L_q(\Omega_0) \times W^1_q(\Omega_0), \quad
\hat \CH_q = \{ (\vf,g) \in \CH_q: \; \int_{\Omega_0} g \,\dx =0 \}, \label{def:Hq}\\
& \CD_q = W^2_q(\Omega_0) \times W^1_q(\Omega_0), \quad
\hat \CD_q= \{(\vf,g) \in \CD_q:\; \int_{\Omega_0}g \,\dx=0 \}. \label{def:Dq}
\end{align}
%
Next, for $0<\ep<\pi$ and $\lambda_0>0$ we introduce 
\begin{align} \label{def:sec}
&\Sigma_\ep = \{ \lambda \in \C \setminus \{0\}: \; |{\rm arg}\lambda|\leq \pi-\ep\}, \quad \Sigma_{\ep,\lambda_0}=\lambda_0+\Sigma_\ep , \quad \Lambda_{\ep,\lambda_0} = \{ \lambda \in \Sigma_\ep: \; |\lambda|\geq\lambda_0 \},\\
&\C_+ := \{ \lambda \in \C: \; {\rm Re}\lambda \geq 0\}.
\end{align}
We also recall the definition of $\CR$-boundedness:
\begin{Definition}\label{dfn:7.1} Let $X$ and $Y$ be two Banach spaces, and 
$\|\cdot\|_X$ and $\|\cdot\|_Y$ their norms.  
A family of operators $\CT \subset \CL(X, Y)$ is called 
$\CR$-bounded on $\CL(X, Y)$ if there exist constants $C > 0$
and $p \in [1, \infty)$ such that 
for any $n \in \BN$, $\{T_j\}_{j=1}^n \subset \CT$ and $\{f_j\}_{j=1}^n
\subset X$, the inequality 
$$\int^1_0\|\sum_{j=1}^n r_j(u)T_jf_j\|_Y^p\,du
\leq C\int^1_0\|\sum_{j=1}^nr_j(u)f_j\|_X^p\,du,
$$
where $r_j: [0, 1] \to \{-1,1\}$, $j \in \BN$, are the Rademacher functions
given by 
$r_j(t) = {\rm sign}(\sin (2^j\pi t))$.  The smallest such
$C$ is called $\CR$-bound of $\CT$ on 
$\CL(X, Y)$ which is written by $\CR_{\CL(X, Y)}\CT$.
\end{Definition}

The Fourier transform and its inverse are defined as  
\begin{equation} \label{def:FT}
\CF[f](\tau)=\int_\R e^{-it\tau}f(t)dt, \quad \CF^{-1}[f](t)=\frac{1}{2\pi}\int_\R e^{it\tau}f(\tau)d\tau   
\end{equation}
We shall also need the Laplace transform and its inverse 
\begin{equation} \label{def:LT}
\CL[f](\lambda)=\int_\R e^{-\lambda t}f(t)dt, \quad \CL^{-1}[f](t)=\frac{1}{2\pi}\int_\R e^{\lambda t}f(\tau)d\ .   
\end{equation}
Finally, by $E(\cdot)$ we shall denote a non-negative non-decreasing continuous function such that $E(0)=0$.  We use the notation $E(T)$ in particular to denote those constants in various inequalities, which can be made arbitrarily small by taking $T$ sufficiently small.
%
\subsection{Main results}

The first main result of this paper gives the local well-posedness for system \eqref{i1a}-\eqref{i1b} with Dirichlet boundary condition. 
\begin{Theorem} \label{t1}
Let $\Omega_0 \subset \R^3$ be a bounded uniform $C^2$ domain. Assume  
$$
\vr_0 \in W^1_q(\Omega_0), \quad  
\vu_0 \in B^{2-2/p}_{q,p}(\Omega_0) 
$$
and 
$$
\vV \in L_p(0,T;W^2_q(\R^3)) \cap W^1_p(0,T;L_q(\R^3))
$$
with $2<p<\infty$, $3<q<\infty$ and $\frac{2}{p}+\frac{3}{q}<1$.
Then for any $L>0$ there exists $T>0$ such that if  
\begin{equation} \label{init:norm}
\|\vr_0\|_{W^1_q(\Omega_0)}+\|\vu_0\|_{B^{2-2/p}_{q,p}(\Omega_0)} + \|\vV\|_{L_p(0,T;W^2_q(\R^3))}+\|\de_t \vV\|_{L_p(0,T;L_q(\R^3))} \leq L
\end{equation}
then the system \eqref{i1a}-\eqref{i7} 
admits a unique strong solution $(\vr,\vu) \in {\cal X}(T)$ and 
\begin{equation} \label{est:t1}
\|\vr,\vu\|_{{\cal X}(T)}\leq CL.    
\end{equation}
\end{Theorem}
\begin{Remark}
Let us comment  the restrictions on $p$ and $q$. The condition $q>3$ is natural
as we shall repeatedly use the embedding $W^1_q(\Omega_0) \subset L_\infty(\Omega_0)$. 
However, a stronger condition $\frac{2}{p}+\frac{3}{q}<1$ is required since we need the embedding $B^{2(1-1/p)}_{q,p}(\Omega_0) \subset W^1_{\infty}(\Omega_0)$ to prove Lemma \ref{l:nonlin2}, see Corollary \ref{c:imbed1}.
\end{Remark}
The second main result gives the global well-posedness:
\begin{Theorem} \label{t2}
Let $\Omega_0 \subset \R^3$ be a bounded uniform $C^2$ domain. Assume that 
$$
\vr_0 \in W^1_q(\Omega_0), \quad  
\vu_0 \in B^{2-2/p}_{q,p}(\Omega_0). 
$$
Furthermore, let $\vr^*,\gamma>0$ be given constants. Then there exists $\epsilon>0$ such that if 
\begin{equation} \label{init:norm2}
\|\vr_0-\vr^*\|_{W^1_p(\Omega_0)}+\|\vu_0-\vV(0)\|_{B^{2-2/p}_{q,p}(\Omega_0)}+ \|e^{\gamma t}(\de_t \vV,\nabla_x \vV, \nabla^2_x \vV)\|_{L_p(0,T;L_q(\R^3))}\leq\epsilon 
\end{equation}
then the unique strong solution to \eqref{i1a}-\eqref{i7} is defined globally in time and 
\begin{align}
&\|e^{\gamma t}\nabla \vr\|_{L_p(0,\infty;L_q(\Omega_t))} + \|e^{\gamma t}\de_t\vr\|_{L_p(0,\infty;W^1_q(\Omega_t))}
+\|\vr-\vr^*\|_{L_\infty(0,\infty;W^1_q(\Omega))} \nonumber\\
&+\|e^{\gamma t}\de_t \vu\|_{L_p(0,\infty;L_q(\Omega_t))}
+\|e^{\gamma t}\nabla_x \vu\|_{L_p(0,\infty;W^1_q(\Omega_t))}
+\|e^{\gamma t}(\vu-\vV)\|_{L_p(0,\infty;L_q(\Omega_t))}
\leq C\epsilon, \label{est:glob1}\\[5pt]     
&\|\vu\|_{L_p(0,\infty;L_q(\Omega_t))}\leq C\epsilon+\|\vV\|_{L_p(0,\infty;L_q(\Omega_t))}. \label{est:glob2}
\end{align}
\end{Theorem}

The paper is structured as follows. In Section 2 we rewrite the problem on a fixed domain using Lagrangian coordinates. { In Section 3 we prove maximal regularity and exponential decay results for linear problems corresponding to linearization of our system in Lagrangian coordinates. As many similar results are already well known, we skip certain technical details referring to relevant papers. In the same section we also recall some imbedding properties.} Section 4 is dedicated to the proof of Theorem \ref{t1}. We reduce the problem to homogeneous boundary condition and show appropriate estimates of the right hand side of the problem in Lagrangian coordinates and conclude using fixed point argument and linear result recalled in Section 3. In Section 5 we prove Theorem \ref{t2}. For this purpose we obtain appropriate estimates of the right hand side which allow to show uniform in time estimate for the solution using exponential decay property of the linear problem. This estimate allow to prolong the solution for arbitrarily large time. 
\section{Lagrangian transformation}
Let us start with a following observation 
\begin{Lemma} \label{l:lag1}
Let $p$ and $q$ satisfy the assumptions of Theorem \ref{t1}. Then\\
(i) if  
$\|f\|_{L_p(0,T;W^2_q(\Omega_0))} \leq M$
for some $M>0$, then
\begin{equation} \label{2:1}
\int_0^T \|\nabla f(t,\cdot)\|_{L_\infty(\Omega_0)} \dt \leq M E(T). \end{equation}
\noindent
(ii) if 
$\|e^{\gamma t}f\|_{L_p(0,\infty;W^2_q(\Omega_0))}\leq M$ 
for some $M,\gamma>0$ then
\begin{equation} \label{2:1a}
\int_0^\infty \|\nabla f(t,\cdot)\|_{L_\infty(\Omega_0)} \dt \leq CM. \end{equation}
\end{Lemma}
\bProof 
By the imbedding theorem and H\"older inequality we have 
\begin{align*}
\int_0^T\|\nabla f(t,\cdot)\|_{L_\infty(\Omega_0)} \dt \leq C \int_0^T \|f(t,\cdot)\|_{W^2_q(\Omega_0)}\dt \leq T^{1/p'}\int_0^T \left(\|f(t,\cdot)\|_{W^2_q(\Omega_0)}^p \right)^{1/p} \dt \leq ME(T),   
\end{align*}
which proves the first assertion, and for the second we have 
\begin{align*}
\int_0^\infty\|\nabla f(t,\cdot)\|_{L_\infty(\Omega_0)} \dt\leq 
C \int_0^\infty e^{-\gamma t}e^{\gamma t}\|f(t,\cdot)\|_{W^2_q(\Omega_0)} \dt\\ \nonumber
\leq \left( \int_0^\infty e^{-\gamma tp'} \dt \right)^{1/p'}
\|e^{\gamma t}f\|_{L_p(0,\infty;W^2_q(\Omega_0))}\leq CM.
\end{align*}
\qed

In order to transform the problem \eqref{i1a}-\eqref{i1b} to a fixed domain we introduce the change of coordinates 
\begin{equation}\label{eq:coc}
\frac{{\rm d}}{{\rm d}t} \vc{X_u}(t, y) = \vc{u} \Big( t, \vc{X_u}(t, y) \Big) \quad \textrm{for}\; t > 0, \quad \vc{X_u}(0, y) = y,
\end{equation}
i.e.
\begin{equation}\label{eq:coc1}
\bX_\vu(t,y)=y+\int_0^t \vu(s,\bXu(s,y)) \d s.
\end{equation}
Then for any differentiable function $f$ defined on $Q_T$ we have 
\begin{equation} \label{dt:lag}
\frac{\d}{\d t} f(t,\bX_{\vu}(t,y))=\frac{\de}{\de t}f(t,\bX_{\vu}(t,y))+\vu\cdot\nabla_x f(t,\bX_\vu(t,y)).
\end{equation}
Let us define transformed density and velocities on a fixed domain $\Omega_0$: 
\begin{equation}
\vrt(t,y)=\vr(t,\bXu(t,y)), \quad \vut(t,y)=\vu(t,\bXu(t,y)), \quad \vVt(t,y) = \vV(t,\bXu(t,y)).  \end{equation}
\begin{Lemma} \label{l:lag2}
Assume that 
\begin{equation} \label{small:2}
\int_0^T \|\nabla_y \vut\|_{\infty} \dt \leq \delta      
\end{equation}
for sufficiently small $\delta>0$. Then the  
inverse to $\bX_{\vu}$, i.e. $\vc{Y}(t,x)$ defined as 
\begin{equation}\label{eq:inv}
\vc{X}_{\vu}(t,\vc{Y}(t,x)) = x \quad \forall t \geq 0, \; x \in \Omega_t,
\end{equation}   
is well defined and its Jacobian can be expressed in a following way  
\begin{equation}\label{2:3}
\nabla_x \bY(t,\bXu(t,y)) = [\nabla_y \bXu(t,y)]^{-1} = \bI + \vE^0(\ku(t,y)),
\end{equation}
where
\begin{equation} \label{ku}
{\bf k}_{\vut}(t,y)=\int_0^t \nabla_y \vut(s,y) \d s
\end{equation}
and $\vE^0(\cdot)$ is a $3\times 3$ matrix of smooth functions with $\vE^0(0) = 0$.
\end{Lemma}
{\bf Proof:} We have 
\begin{equation}\label{2:2}
\frac{\pd X_i}{\pd y_j}(t,y) = \delta_{ij} 
+ \int^t_0\frac{\pd \tilde{u}_i}{\pd y_j}(s,y)\,ds.
\end{equation}
Therefore, if $\eqref{small:2}$ holds for sufficiently small $\delta$ then 
$\bY(t,x)$ is well defined and we have \eqref{2:3}-\eqref{ku} with $\vE^0$ 
as in the statement of the Lemma. Next, by the boundary condition \eqref{i6a} 
we have 
$$
\bX_{\vu}(\Gamma_0,t)=\Gamma_t \quad {\rm for} \; t>0
$$
and
$$
\bX_{\vu}(y,t) \subset \Omega_t \quad {\rm for} \; t>0,\, y \in \Omega_0.  
$$
Finally, it is well known that $\bX_\vu$ is a diffeomorphism which completes the proof.

\qed

Note that by \eqref{2:3} we can write
\begin{equation} \label{dx:lag}
\nabla_x = [\bI + \vE^0(\ku)]\nabla_y.
\end{equation}
\begin{Lemma}
Let $(\vr,\vu)$ be a solution to \eqref{i1a}-\eqref{i7}. Then $(\vrt,\vut)$ solve the following system of equations on the fixed domain $\Omega_0$
\begin{align} \label{ME:lag1}
&\vrt \vut_t - \mu\Delta_y\vut-\left(\frac{\mu}{3}+\zeta \right)\nabla_y\div_y\vut+\nabla_y \pi(\vrt)= \vF(\vrt,\vut),\\
\label{CE:lag1}
&\vrt_t + \vrt \div_y \vut = G(\vrt,\vut),\\
&\label{icbc:lag1} \vut|_{t=0}=\vu_0, \quad \vut|_{\de \Omega_0}=\vVt.
\end{align}
The $i$-th component of $\vF(\cdot,\cdot)$ is given by
\begin{align} \label{F:lag1}
&F_i(\vrt,\vut)=-E^0_{ij}\de_{y_j}\pi(\vrt) + R_i(\vut), 
\end{align}
and
\begin{equation}
\label{G:lag1}
G(\vrt,\vut)= -\vrt \vE^0_{ij}(\ku)\frac{\de \tilde u_i}{\de y_j},
\end{equation}
where the components $R_i(\cdot)$ of $\vR(\cdot)$ are expressed as
\begin{equation} \label{def:R}
R_i(\vut)=\mu[ A_{2\Delta}(\ku)\nabla^2_y\vut
+ A_{1\Delta}(\ku)\nabla_y\vut]_i 
+\left(\frac{\mu}{3}+\zeta \right)[A_{2\div, i}(\ku)\nabla^2_y\vut
+ A_{1\div, i}(\ku)\nabla_y\vut],
\end{equation}
with $A_{j\Delta}$ and $A_{j\div}$ ($j=1,2$) given in \eqref{a2delta}, \eqref{a1delta}, \eqref{a2div} and \eqref{a1div}, respectively. 
\end{Lemma}
\bProof
We have  
\begin{equation}\label{lag:div}
\div_x \vu = \div_y \vut + \vE^0:\nabla_y \vut,
\end{equation}
where $\vE^0:\nabla_y \vut = \vE^0_{ij}(\ku)\frac{\de \tilde u_i}{\de y_j}$,
which together with \eqref{dt:lag} gives \eqref{CE:lag1}.

In order to transform the momentum equation \eqref{i1b}
it is convenient to rewrite it, using \eqref{i1a} and \eqref{i4}, as 
\begin{equation} \label{i1bb}
\vr (\de_t \vu + \vu\cdot \nabla_x \vu)-\mu\Delta_x \vu-(\frac{\mu}{3}+\zeta)\nabla\div_x \vu + \nabla_x \pi(\vr)=0.
\end{equation}
We have
\begin{align} \label{lag:5}
&\de_{x_i}\pi(\vr)
= \de_{y_i}\pi(\vrt)+E^0_{ij}\de_{y_j}\pi(\vrt).
\end{align}
Now we need to transform second order operators.   
By \eqref{dx:lag}, we have
$$
\Delta_x \vu = \frac{\pd}{\pd x_k}\lr{\frac{\pd \vu}{\pd x_k}}
= \lr{\delta_{kl} + \vE^0_{kl}(\ku)}
\frac{\pd}{\pd y_l}
\lr{\lr{\delta_{km} + \vE^0_{km}(\ku)}\frac{\pd \vut}{\pd y_m}}.
$$
Therefore 
\begin{equation} \label{lag:6}
\Delta_x \vu = \Delta_y \vut + A_{2\Delta}(\ku)\nabla^2_y\vut
+ A_{1\Delta}(\ku)\nabla_y\vut
\end{equation}
with
\eq{
A_{2\Delta}(\ku)\nabla^2_y\vut &= 2 \sum_{l,m} \vE^0_{kl}(\ku)
\frac{\pd^2\vut}{\pd y_l\pd y_m}
+ \sum_{k,l,m} \vE^0_{kl}(\ku)\vE^0_{km}(\ku)
\frac{\pd^2\vut}{\pd y_l \pd y_m}, \label{a2delta} }
\eq{
A_{1\Delta}(\ku)\nabla_y\vut  = &(\nabla_{\ku} \vE^0_{l m})(\ku)
\int^t_0(\pd_l\nabla_y\vut)\,\d s \frac{\pd \vut}{\pd y_m}\\
&+ \vE^0_{kl}(\ku) (\nabla_{\ku} \vE^0_{km})(\ku)
\int^t_0\pd_l\nabla_y\vut\, \d s\frac{\pd \vut}{\pd y_m}. \label{a1delta}
}
%
Next, by \eqref{lag:div}
$$
\frac{\pd}{\pd x_i}\div_x\vu 
= \sum_{k=1}^3(\delta_{ik} + \vE^0_{ik}(\ku))\frac{\pd}{\pd y_k}
\lr{\dv_y\vut + \sum_{l, m=1}^3 \vE^0_{l m}(\ku)\frac{\pd \tilde u_l}{\pd y_m}},
$$
so we obtain
\begin{equation} \label{lag:7}
\frac{\pd}{\pd x_i}\dv_x\vu
= \frac{\pd}{\pd y_i}\dv_y\vut + A_{2\div, i}(\ku)\nabla^2_y\vut
+ A_{1\div, i}(\ku)\nabla_y\vut,
\end{equation}
where
\eq{
A_{2\div,i}(\ku)\nabla^2_y\vut
& = \sum_{l, m=1}^3 \vE^0_{l m}(\ku)\frac{\pd^2 \tilde u_l}{\pd y_m \de y_i}
+ \sum_{k=1}^3 \vE^0_{ik}(\ku)\frac{\pd}{\pd y_k}\dv_y\vut
+ \sum_{k, l=1}^3 \vE^0_{ik}(\ku)\vE^0_{l m}(\ku)
\frac{\pd^2 \tilde u_l}{\pd y_k \pd y_m}, \label{a2div}
}
\eq{
A_{1\div, i}(\ku)\nabla_y\vut
 =& \sum_{l, m=1}^3(\nabla_{\ku} \vE^0_{l m})(\ku)
\int^t_0\pd_i\nabla_y\vut\,\d s\frac{\pd \tilde u_l}{\pd y_m} \\
&+ \sum_{k,l, m=1}^3 \vE^0_{ik}(\ku)(\nabla_{\ku} \vE^0_{l m})(\ku)
\int^t_0\pd_k\nabla_y\vut\,\d s\frac{\pd \tilde u_l}{\pd y_m}. \label{a1div}
}

Putting together \eqref{lag:5}, \eqref{lag:6} and \eqref{lag:7} gives \eqref{ME:lag1} with \eqref{F:lag1}.

\qed

\section{Linear theory and auxiliary results}
\subsection{$L_p-L_q$ maximal regularity}
This section is dedicated to $L_p-L_q$ maximal regularity results for linear problems on a fixed domain, which will be used in the proof of Theorem \ref{t1}. 
We shall rely on similar already known results, therefore we will be able to avoid most technicalities referring to relevant works and presenting only an outline of the proof. The linearized system of equations on the fixed domain $\Omega_0$ reads as
\begin{align} 
&\vr_0\vu_t - \mu \Delta_y \vu - (\frac{\mu}{3}+\zeta) \nabla_y \div_y \vu + \gamma_0 \nabla_y \eta = \vf, \label{ME:lin0} \\
&\eta_t + \vr_0 \div_y \vu = g, \label{CE:lin0}\\
&\vu|_{\de \Omega_0}=0, \quad \vu|_{t=0}=\vu_0, \label{icbc:lin0}
\end{align}
where the unknowns are $\eta$ and $\vu$, while $\vr_0, \gamma_0, \vf, g$ are given functions such that $\vr_0 \geq c>0$ and $\gamma_0\geq 0$.
To show the local well-posedness for the Dirichlet boundary condition we will use the following result
\begin{Theorem} \label{p1}
Let $1<p,q<\infty$ and $\frac{2}{p}+\frac{1}{q} \neq 1$.
Let $\vr_0,\vu_0,\mu$ and $\zeta$ satisfy the assumptions of Theorem \ref{t1}. 
Moreover, let $\Omega_0 \subset \R^n$ be a uniform $C^2$ domain.
If $\frac{2}{p}+\frac{1}{q}<1$, assume additionally that the initial velocity satisfies the compatibility condition $\vu_0|_{\partial\Omega_0}=0$.
Finally, assume that for some $T>0$ 
$$
\vf \in L_p(0,T;L_q(\Omega_0)),\quad g \in L_p(0,T;W^1_q(\Omega_0)). 
$$
Then the problem \eqref{ME:lin0}-\eqref{icbc:lin0} admits a unique solution $(\eta,\vu) \in {\cal Y}(T)$ such that 
\begin{align} \label{est:linME1_A}
\|\eta,\vu\|_{{\cal Y}(T)}\leq C(T,\mu,\zeta,\|\vr_0\|_{L_{\infty}(\Omega_0)}) e^{\gamma_0 t}[\|\vu_0\|_{B^{2-2/p}_{q,p}(\Omega_0)}+\|\vf\|_{L_p(0,T;L_q(\Omega_0)}+\|g\|_{L_p(0,T;W^1_q(\Omega_0)}]
\end{align} 
for some positive constant $\gamma_0$.
\end{Theorem}  
In order to remove inhomogeneity from the boundary condition we also need to consider separately the linearized momentum equation 
\begin{align} 
&\vr_0\vu_t - \mu \Delta_y \vu - (\frac{\mu}{3}+\zeta) \nabla_y \div_y \vu = \vf, \label{ME:lin} \\
&\vu|_{\de \Omega_0}=0, \quad \vu|_{t=0}=\vu_0. \label{icbc:lin}
\end{align}
For the latter we have 
\begin{Theorem} \label{p1b}
Let $\Omega_0,\vu_0,p,q,\mu,\zeta$ and $\vf$ satisfy the assumptions of Theorem \ref{p1}. Then the problem \eqref{ME:lin}-\eqref{icbc:lin} admits a unique solution $\vu \in L_p(0,T;W^2_q(\Omega_0)) \cap W^1_p(0,T;L_q(\Omega_0))$ such that 
\begin{align} \label{est:linME1_B}
\|\vu\|_{L_p(0,T;W^2_q(\Omega_0)}+\|\vu_t\|_{L_p(0,T;L_q(\Omega_0)} \leq C(T,\mu,\zeta,\|\vr_0\|_{L_{\infty}(\Omega_0)}) [\|\vu_0\|_{B^{2-2/p}_{q,p}(\Omega_0)}+\|\vf\|_{L_p(0,T;L_q(\Omega_0)}].
\end{align} 
\end{Theorem}  

The proofs of Theorems \ref{p1} and \ref{p2} are based on the concept of $\CR$ - boundedness. For a space $X$, by $\CD(\R,X)$ we shall denote $C^\infty$ functions with compact support and values in $X$ and by $\CS(\R,X)$ a Schwartz space of $X$-valued functions.  
A fundamental tool enabling application of $\CR$-boundedness to maximal regularity of PDE is the following Weis' vector-valued Fourier multiplier theorem (\cite{Weis}).
\begin{Theorem} \label{thm:Weis}
Let $X$ and $Y$ be UMD spaces and $1<p<\infty$. Let $M \in C^1(\R\setminus \{0\}, \CL(X,Y))$. Let us define the operator 
$T_M:\CF^{-1} \CD(\R,X) \to \CS'(\R,Y)$:
\begin{equation} \label{def:TM}
T_M\phi(\tau) = \CF^{-1}[M\CF[\phi](\tau)].    
\end{equation}
Assume that 
\begin{equation} \label{rbound:Weis}
\CR_{\CL(X,Y)}(\{M(\tau): \; \tau \in \R \setminus \{0\}\})=\kappa_0<\infty, \quad 
\CR_{\CL(X,Y)}(\{\tau M'(\tau): \; \tau \in \R \setminus \{0\}\})=\kappa_1<\infty.
\end{equation}
Then, the operator $T_M$ defined in \eqref{def:TM} is extended to a bounded linear operator $L_p(\R,X) \rightarrow L_p(\R,X)$ and 
$$
\|T_M\|_{\CL(L_p(\R,X),L_p(\R,Y))} \leq C(\kappa_0+\kappa_1),
$$
where $C=C(p,X,Y)>0$. 
\end{Theorem}
\begin{Remark}
For definitions and properties of UMD spaces we refer the reader for example to Chapter 4 in \cite{HNVW}. Here let us only note that $L_p$ spaces and $W^k_p$ spaces are UMD for $1<p<\infty$.
\end{Remark}
The following result (Theorem 2.17 in \cite{ES}) explains how to obtain $L_p$ - maximal regularity using Theorem \ref{thm:Weis} and Laplace transform:  
\begin{Theorem} \label{thm:maxreg}
Let $X$ and $Y$ be UMD Banach spaces and $1<p<\infty$. Let $0<\ep<\frac{\pi}{2}$ and $\gamma_1 \in \R$. 
Let $\Phi_\lambda$ be a $C^1$ function of $\tau \in \R\setminus\{0\}$ where $\lambda=\gamma+i\tau \in \Sigma_{\ep,\gamma_1}$ with values in $\CL(X,Y)$. 
Assume that 
\begin{equation}  
\CR_{\CL(X,Y)}(\{\Phi_\lambda: \; \lambda \in \Sigma_{\ep,\gamma_1}\})\leq M,\quad
\CR_{\CL(X,Y)}\left(\left\{ \tau\frac{\de}{\de\tau}\Phi_\lambda: \; \lambda \in \Sigma_{\ep,\gamma_1}\right\}\right) \leq M
\end{equation}
for some $M>0$. Let us define 
\begin{equation} \label{def:Psi}
\Psi f(t)=\CL^{-1}[\Phi_\lambda \CL[f](\lambda)](t) \quad {\rm for} \quad f \in \CF^{-1}\CD(\R,X),
\end{equation}
where $\CL$ and $\CL^{-1}$ are the Laplace transform and its inverse defined in
\eqref{def:LT}.
Then 
\begin{equation} \label{est:maxreg}
\|e^{-\gamma t}\Psi f\|_{L_p(\R,Y)} \leq C(p,X,Y)M\|e^{-\gamma t}\|_{L_p(\R,X)} \quad \forall \gamma \geq \gamma_1  
\end{equation}
\end{Theorem}
{\bf Proof:} For $\lambda=\gamma+i\tau$ we have the following relation between 
Laplace and Fourier transforms defined in, respectively, \eqref{def:LT} and \eqref{def:FT}: 
\begin{align*}
&\CL[f](\lambda)=\int_\R e^{-\lambda t}f(t)dt = \CF[e^{-\gamma t}f](\tau),\\
&\CL^{-1}[f](t)=\frac{1}{2\pi}\int_\R e^{\lambda t}f(\lambda)d\tau=e^{\gamma t}\CF^{-1}[f](t)
\end{align*}
Therefore by \eqref{def:Psi} we have
$$
e^{-\gamma t}\Psi f(t)=\CF^{-1}[\Psi_{\gamma+i\tau}\CF[e^{-\gamma t}](\tau)](t).
$$
Applying Theorem \ref{thm:Weis} to the above formula we conclude \eqref{est:maxreg}

\qed

If $(\vu,\eta)$ solves \eqref{ME:lin0}-\eqref{icbc:lin0} then 
$\vu_\lambda=\CL[\vu](\lambda)$, $\eta_\lambda=\CL[\eta](\lambda)$ satisfy the following resolvent problem 
\begin{align} 
&\vr_0\lambda \vu_\lambda - \mu \Delta_y \vu_\lambda - (\frac{\mu}{3}+\zeta) \nabla_y \div_y \vu_\lambda + \gamma_0 \nabla_y \eta_\lambda = \hat \vf, \label{ME:res0} \\
&\lambda \eta_\lambda + \vr_0 \div_y \vu_\lambda = \hat g, \label{CE:res0}\\
&\vu_\lambda|_{\de \Omega_0}=0, \label{bc:res0}  
\end{align}
with $\hat \vf = \CL[\vf](\lambda)$ and $\hat g = \CL[g](\lambda)$. Therefore if we define by $\Phi_\lambda$ a solution operator to \eqref{ME:res0}-\eqref{bc:res0}, i.e. 
$$
(\vu_{\lambda},\eta_\lambda) = \Phi_\lambda(\tilde \vf,\tilde g)
$$
then
$(\vu,\eta)=\Psi (\vf,g)$, where $\Psi$ is defined in \eqref{def:Psi} satisfies the original problem \eqref{ME:lin0}-\eqref{icbc:lin0}.
Hence in order to prove Theorem \ref{p1} it suffices to show ${\cal R}$-boundedness of a family of solutions to \eqref{ME:res0}-\eqref{CE:res0}. 
The result is

\begin{Theorem} \label{thm:Rbound}
Let $1 < p < \infty$ and $0 < \epsilon < \pi/2$.
Assume that $\Omega_0$ is a uniform $C^2$ domain. 
Then, there exist a positive constant $\lambda_0$ and 
operator families $\CA(\lambda) \in {\rm Hol}\,(\Lambda_{\epsilon, \lambda_0},
\CL(L_q(\Omega_0)\times W^1_q(\Omega_0), W^1_q(\Omega_0)))$ and  
$\CB(\lambda) \in {\rm Hol}\,(\Lambda_{\epsilon, \lambda_0},
\CL(L_q(\Omega_0)\times W^1_q(\Omega_0), W^2_q(\Omega_0)^N))$, 
such that for any $(\hat \vf, \hat g) \in L_q(\Omega_0) \times W^1_q(\Omega_0)$
and $\lambda \in \Lambda_{\epsilon, \lambda_0}$, 
$(\eta_\lambda = \CA(\lambda)(\hat \vf,\hat g), 
\vu_\lambda= \CB(\lambda)(\hat \vf, \hat g))$
is a unique solution of \eqref{ME:res0}-\eqref{bc:res0} and 
\begin{align*}
\CR_{\CL(L_q(\Omega_0) \times W^1_q(\Omega_0), W^1_q(\Omega_0))}(\{(\tau\pd_\tau)^\ell
\CA(\lambda):\; \lambda \in \Sigma_{\epsilon, \lambda_0}\}) &\leq M, \\
\CR_{\CL(L_q(\Omega_0) \times W^1_q(\Omega_0), H^{2-j}_q(\Omega_0)^N)}(\{(\tau\pd_\tau)^\ell
(\lambda^{j/2}\CB(\lambda)):\; 
\lambda \in \Sigma_{\epsilon, \lambda_0}\}) &
\leq M
\end{align*}
for $\ell=0,1$, $j=0,1,2$ and some constant $M>0$, where { by ${\rm Hol}$ we denote the space of holomorphic operators.}
\end{Theorem}
{
\begin{Remark}
We say that the operator valued function $T(\lambda)$ is holomorphic if it is differentiable in norm for all $\lambda$ in a complex domain. For more details, see \cite[Chaper VII,\$1.1]{Kato}.
\end{Remark}
}
A particular advantage of this approach is a direct treatment of the continuity equation. Namely, in order to prove Theorem \ref{thm:Rbound} we use \eqref{CE:res0} to get $\eta_\lambda = \lambda^{-1}(g-\vr_0 \div_y \vu_\lambda)$. Plugging the latter to \eqref{ME:res0} and skipping the subscripts $\lambda$ we get 
\begin{equation} \label{ME:res0b}
\vr_0\lambda \vu - \mu \Delta_y \vu - (\frac{\mu}{3}+\zeta+\gamma_0\lambda^{-1}\vr_0) \nabla_y \div_y \vu  = \vf - \gamma_0 \lambda^{-1}(\nabla g + \nabla \vr_0 \div\vu) , \quad
\vu|_{\de \Omega_0}=0    
\end{equation}
and the problem is reduced to showing ${\cal R}$-boundedness for the family of solutions to \eqref{ME:res0b}. With this approach Proposition \ref{p1b} becomes a special case of Theorem \ref{p1}. Namely, resolvent problem corresponding to \eqref{ME:lin} reads
\begin{equation} \label{ME:res}
\vr_0\lambda \vu - \mu \Delta_y \vu - (\frac{\mu}{3}+\zeta) \nabla_y \div_y \vu = \vf, \quad
\vu|_{\de \Omega_0}=0.    
\end{equation}
As the second term on the right hand side of \eqref{ME:res0b}$_1$ is of lower order, we can 
reduce \eqref{ME:res0b} to 
\begin{equation} \label{ME:res0c}
\vr_0\lambda \vu - \mu \Delta_y \vu - (\frac{\mu}{3}+\zeta+\gamma_0\lambda^{-1}\vr_0) \nabla_y \div_y \vu  = \vF , \quad
\vu|_{\de \Omega_0}=0.    
\end{equation}
Now \eqref{ME:res} is a particular case of \eqref{ME:res0c} with $\gamma_0 = 0$.
The following result gives ${\cal R}$-boundedness for \eqref{ME:res0c}: 

\begin{Proposition} \label{prop:Rbound} Let $1 < q < \infty$ and 
$0 < \epsilon < \pi/2$.  Assume that $\Omega_0$ is a uniform 
$C^2$ domain in $\R^N$.  Then, there exists a positive constant
$\lambda_0$ such that there exists   
an operator family $\CC(\lambda) \in {\rm Hol}\,
(\Sigma_{\epsilon, \lambda_0}, \CL(L_q(\Omega_0)^N, H^2_q(\Omega_0)^N))
$
such that for any 
$\lambda \in \Lambda_{\epsilon, \lambda_0}$
and $\vF \in L_q(\Omega_0)^N$, 
$\bv = \CC(\lambda)\vF$ is a unique 
solution of \eqref{ME:res0c}, and
$$\CR_{\CL(L_q(\Omega_0)^N, H^{2-j}_q(\Omega_0)^N)}
(\{(\tau\pd_\tau)^\ell\CC(\lambda) \mid \lambda \in 
\Sigma_{\epsilon, \lambda_0}\}) \leq M$$
for $\ell=0,1$, $j=0,1,2$ and some constant $M>0$.
\end{Proposition}
{\bf Proof.} An analog of Proposition \ref{prop:Rbound} has been shown in (\cite{ES}, Theorem 2.10)  for a problem 
\begin{equation} \label{ME:resES}
\lambda \vu - \mu \Delta_y \vu - (\nu+\gamma^2\lambda^{-1}\vr_0) \nabla_y \div_y \vu  = \vF , \quad
\vu|_{\de \Omega_0}=0    
\end{equation}
where $\mu,\nu$ and $\gamma$ are constants satisfying $\mu+\nu>0$ and 
$\gamma>0$. The proof requires only minor modifications in order to prove Proposition \ref{prop:Rbound}, therefore we present only a sketch. First we solve a problem with constant coefficients in the whole space
\begin{equation} \label{ME:res:const}
\vr_0^*\lambda \vu - \mu \Delta_y \vu - (\frac{\mu}{3}+\zeta+\gamma_0^*\lambda^{-1}\vr_0^*) \nabla_y \div_y \vu  = \vF \quad \in \R^n,
\end{equation}
where $\rho_0^*>0$ and $\gamma_0^* \geq 0$ are constants. $\CR$-boundedness
for \eqref{ME:res:const} can be shown following the proof Theorem 3.2 in \cite{ES}, the only difference is that now we have $\gamma_0^* \rho_0^*$ instead of $\gamma^2$. However, strict positivity of this constant neither its square structure is not necessary, nonnegativity is sufficient. 

Next we consider \eqref{ME:res:const} in a half-space supplied with the boundary condition $u|_{\de \R^n_+}=0$. Here we can follow the proof of Theorem 4.1 in \cite{ES} which works without modifications for $\gamma_0^* \rho_0^* \geq 0$. 

The third step consists in showing $\CR$-boundedness in a bent half-space, the necessary result is Theorem 5.1 in \cite{ES}, where replacing $\gamma^2>0$ with  
$\gamma_0^* \rho_0^* \geq 0$ is again harmless. 

The final step is an introduction of a partition of unity and application of properties of a uniform $C^2$ domain. Here we have to deal with variable coefficients which is not in the scope of Theorem 2.10 in \cite{ES}. 
However, we can refer to a more recent result, Theorem 4.1 in \cite{PSZ} which gives $\CR$-boundedness for a resolvent problem corresponding to more complicated system describing flow of a two-component mixture.
It is enough to follow step by step Section 6.3 of \cite{PSZ} omitting all terms which are not relevant here. 

\qed
\noindent
{\bf Proof of Theorem \ref{thm:Rbound}}.
By definition of $\Lambda_{\ep,\lambda_0}$ we have $|\lambda^{-1}|\leq C$ for $\lambda \in \Lambda_{\ep,\lambda_0}$. Therefore 
\begin{align*}
&\|\vf-\gamma_0\lambda^{-1}(\nabla g+\nabla \vr_0\div \vu)\|_{L_p(\Omega_0)}
\leq C[\|f\|_{L_q(\Omega_0)}+\|\nabla \vr_0\|_{L_q(\Omega_0)}\|\div \vu\|_{L_\infty(\Omega_0)})] \leq C [\|f\|_{L_q(\Omega_0)}+\|\vu\|_{W^2_q(\Omega_0)}],\\
& \|\lambda^{-1}(g-\nabla \vr_0\div\vu)\|_{W^1_q}\leq C[\|g\|_{W^1_q(\Omega_0)}+\|\vu\|_{W^2_q(\Omega_0)}],
\end{align*}
and so Theorem \ref{thm:Rbound} results directly from Proposition \ref{prop:Rbound}.

\qed
\smallskip

\noindent
Now it is straightforward to conclude main results of this section.
\smallskip

\noindent 
{\bf Proof of Theorems \ref{p1} and \ref{p2}}.

\noindent Theorem \ref{p1} follows from Theorem \ref{thm:Rbound} and Theorem \ref{thm:maxreg} with $X=L_q(\Omega_0) \times W^1_q(\Omega_0)$ and $Y=W^1_q(\Omega_0)\times W^2_q(\Omega_0)$. Theorem \ref{p2} follows from Proposition \ref{prop:Rbound} and Theorem \ref{thm:maxreg}, this time with $X=L_q(\Omega_0)$ and $Y=W^2_q(\Omega_0)$. 

\qed

\subsection{Exponential decay}

\subsubsection{Linearization}

In order to show the global well-posedness in Theorem \ref{t2} 
we will linearize the problem around the constant $\vr^*$, therefore we consider on the fixed domain $\Omega_0$ a linear problem
\begin{align} 
&\label{ME:lin1a}\vr^*\vu_t - \mu \Delta_y \vu - (\frac{\mu}{3}+\zeta) \nabla_y \div_y \vu +\gamma_0^* \nabla_y \eta= \vf,   \\
&\eta_t + \vr^* \div_y \vu = g, \label{CE:lin1a}\\
&\vu|_{\de \Omega_0}=0, \quad \vu|_{t=0}=\vu_0. \label{icbc:lin1a}
\end{align}
Here again $\eta$ and $\vu$ are unknowns and $f,g$ are given functions, but this time $\vr^*>0$ and $\gamma_0^* \geq 0$ are constants.
The main result of this section is the following exponential decay estimate
\begin{Proposition} \label{p2}
Assume $\vr^*>0$ is a constant and $\Omega_0,p,q,\mu,\zeta,\vu_0$ satisfy the assumptions of Theorem \ref{p1}. Assume moreover that there exists $\gamma>0$ such that
$$
e^{\gamma t} \vf \in L_p(0,\infty;L_q(\Omega_0)), \quad 
e^{\gamma t} g \in L_p(0,\infty;W^1_q(\Omega_0)).
$$
Then \eqref{ME:lin1a}-\eqref{icbc:lin1a} admits a unique solution $\vr,\vu$ such that
\begin{align} \label{est:decay1}
&\|e^{\gamma t}\vu_t\|_{L_p(0,\infty;L_q(\Omega_0))}
+\|e^{\gamma t}\vu\|_{L_p(0,\infty;W^2_q(\Omega_0))}
+\|e^{\gamma t}\nabla \eta\|_{L_p(0,\infty;L_q(\Omega_0))}
+\|e^{\gamma t}\de_t \eta\|_{L_p(0,\infty;W^1_q(\Omega_0))}
\nonumber\\
&\leq C_{p,q}\left(\|\vu_0\|_{B^{2-2/p}_{q,p}(\Omega_0)}+\|e^{\gamma t}\vf\|_{L_p(0,\infty;L_q(\Omega_0))}+\|e^{\gamma t}g\|_{L_p(0,\infty;W^1_q(\Omega_0))}\right).
\end{align}
\end{Proposition} 
As before, we also need a decay estimate for the linear momentum equation 
\begin{align} 
&\label{ME:lin4}\vr^*\vu_t - \mu \Delta_y \vu - (\frac{\mu}{3}+\zeta) \nabla_y \div_y \vu 
= \vf,   \\
&\vu|_{\de \Omega_0}=0, \quad \vu|_{t=0}=\vu_0. 
\end{align}
\begin{Proposition} \label{p2b}
Assume $\vr^*>0$ is a constant and $\Omega_0,p,q,\mu,\zeta,\vu_0$ satisfy the assumptions of Theorem \ref{p1} and $\vf$ satisfies 
the assumptions of Proposition \ref{p2}.
Then \eqref{ME:lin4} admits a unique solution $\vu$ such that
\begin{equation}
\|e^{\gamma t}\vu_t\|_{L_p(0,\infty;L_q(\Omega_0))}
+\|e^{\gamma t}\vu\|_{L_p(0,\infty;W^2_q(\Omega_0))}
\leq C_{p,q}\left(\|\vu_0\|_{B^{2-2/p}_{q,p}(\Omega_0)}+\|e^{\gamma t}\vf\|_{L_p(0,\infty;L_q(\Omega_0))}\right).
\end{equation}
\end{Proposition}

\subsubsection{Proof of Propositions \ref{p2} and \ref{p2b}}

Let us focus on the proof of Proposition \ref{p2}, as the proof Proposition \ref{p2b} is just a simplification of the latter.  Roughly speaking, the point is that the resolvent estimate from the previous section holds for $\lambda$ bounded away from zero and now we have to extend it to cover some range of negative values of ${\rm Re}\, \lambda$ under stronger assumptions (fixed coefficients).
The resolvent problem corresponding to \eqref{ME:lin1a}-\eqref{icbc:lin1a} reads 
\begin{align} 
&\label{res:ME:lin1a} \vr^*\lambda \vu - \mu \Delta_y \vu - (\frac{\mu}{3}+\zeta) \nabla_y \div_y \vu + \gamma_0^* \nabla_y \eta = \vf,   \\
&\lambda \eta + \vr^* \div_y \vu = g, \label{res:CE:lin1a}\\
&\vu|_{\de \Omega_0}=0. \label{res:bc}
\end{align}
\begin{Proposition} \label{p5}
Let $q,\vr_{*},\gamma_{0*},\mu,\zeta$ satisfy the assumptions of Theorem \ref{p1}.
Then there exists $\delta,\ep_1>0$ such that for any 
$$
\lambda \in \Sigma_{\ep_1,-\delta} 
$$
and $(\vf,g) \in \hat \CH_q$ the problem \eqref{res:ME:lin1a}-\eqref{res:bc} admits a unique solution $(\vu,\eta)$ with the estimate
\begin{equation} \label{est:res}
(|\lambda|+1)\|(\vu,\eta)\|_{\CH_q}+\|\vu\|_{W^2_q(\Omega_0)}\leq C \|(\vf,g)\|_{\CH_q}    
\end{equation}
\end{Proposition}
{\bf Proof. Step 1: $\lambda \in \C_+ \cup \Lambda_{\ep,\lambda_0}$.}  The resolvent problem \eqref{res:ME:lin1a}-\eqref{res:bc} is a particular case of
\eqref{ME:res0}-\eqref{bc:res0}
with fixed coefficients. Therefore from Theorem \ref{thm:Rbound} it follows that there exists $\lambda_0 > 0$ such that the solution with the estimate \eqref{est:res} exists for $\lambda \in \Lambda_{\epsilon,\lambda_0}$. It remains to prove the unique existence and \eqref{est:res} for $\{\lambda \in \C_+ : \; |\lambda|\leq \lambda_0\}$. 
Consider first $\lambda \neq 0$. Computing $\eta$ from \eqref{res:CE:lin1a} and plugging to \eqref{res:ME:lin1a} we obtain 
\begin{equation} \label{res:3}
\lambda \vu - \vr_*^{-1}[\mu \Delta\vu+(\frac{\mu}{3}+\zeta+\lambda^{-1}\vr_*)\nabla\div\vu] = \vr_*^{-1}(f-\lambda^{-1}\nabla g)=: \tilde \vf, \quad \vu|_{\de \Omega_0}=0.    
\end{equation}
Let us denote 
\begin{equation} \label{def:alambda}
A_{\lambda}\vu=\vr_*^{-1}[\mu \Delta\vu+(\frac{\mu}{3}+\zeta+\lambda^{-1}\vr_*)\nabla\div\vu]    
\end{equation}
and consider an auxiliary problem with a resolvent parameter $\tau>0$:
\begin{equation} \label{res:tau}
\tau \vv_\tau - A_\lambda \vv_\tau = \tilde \vf.    
\end{equation}
By (\cite{ST}, Theorem 7.1) there exists $\tau_0>0$ such that for $\tau \geq \tau_0$ there exists $(\tau \I - A_\lambda)^{-1} \in \CL(L_q(\Omega_0),\CD_q)$ and 
$$
\tau \|\vv_\tau\|_{L_q(\Omega_0)}+\tau^{1/2}\|\vv_\tau\|_{W^1_q(\Omega_0)}+\|\vv_\tau\|_{W^2_q(\Omega_0)}\leq C \|\tilde \vf\|_{L_q(\Omega_0)}.
$$
Now we rewrite \eqref{res:3} as 
$$
(\lambda-\tau)\vu+(\tau \I-A_\lambda)\vu=\tilde \vf .
$$
Applying $(\tau \I-A_\lambda)^{-1}$ to the above equation we get 
\begin{equation} \label{res:4}
\vu+(\lambda-\tau)(\tau \I-A_\lambda)^{-1}\vu=(\tau\I-A_\lambda)^{-1} \tilde \vf .
\end{equation}
Since $(\lambda-\tau)(\tau \I-A_\lambda)^{-1}$ is compact on $L_q(\Omega_0)$,
in order to prove unique solvability of \eqref{res:4} it is enough to show 
that 
\begin{equation} \label{ker}
{\rm Ker}(\I+(\lambda-\tau)(\tau\I-A_\lambda)^{-1})=\{0\}.
\end{equation}
Assume that 
$$
\left[\I+(\lambda-\tau)(\tau\I-A_\lambda)^{-1}\right] {\bf g}=0.
$$
Setting $\vv^0=(\tau\I-A_\lambda)^{-1}{\bf g}$ we get
$$
(\lambda\I - A_\lambda)\vv^0=0, 
$$
Under the assumptions on the coefficients, by definition of $A_\lambda$ easily verify that $\vv^0=0$, therefore $\bf g=0$ and so \eqref{ker} holds. 
This completes the proof of the unique solvability of \eqref{res:ME:lin1a}-\eqref{res:bc} for $\lambda\in \C_+ \setminus \{0\}.$ 

For $\lambda=0$ problem \eqref{res:ME:lin1a}-\eqref{res:bc} reduces to inhomogeneous Stokes problem 
\begin{align} 
&- \mu \Delta_y \vu + \gamma_0^* \nabla_y \eta = \vf + (\frac{\mu}{3}+\zeta)\vr_*^{-1} \nabla_y g, \label{stokes:1}\\ 
&\div_y \vu =  g, \quad \vu|_{\de \Omega_0}=0.
\label{stokes:2}
\end{align}
In \cite{FS} the unique existence was shown for a problem 
$$
\kappa \vu - \mu \Delta_y \vu + \gamma_0^* \nabla_y \eta = \vf, \quad \div_y \vu =  g, \quad \vu|_{\de \Omega_0}=0
$$
with $(\vf,g) \in \hat \CH_q$ for sufficienlty large $\kappa$. 
By Fredholm alternative, this result together with uniqueness for \eqref{stokes:1}-\eqref{stokes:2}, which is obvious, implies the unique 
solvability of \eqref{stokes:1}-\eqref{stokes:2}. Notice that here we require that $g$ has zero mean. This completes the proof for $\lambda \in \C_+ \cup \Lambda_{\ep,\lambda_0}$ 

{\bf Step 2.} Proof for $\lambda \in \Sigma_{\ep_1,-\delta}$. The crucial observation is that the resolvent set is open. Therefore, since the set $\{\lambda \in \C: {\rm Re} \lambda = 0, \; {\rm Im}\lambda \leq \lambda_0 \}$ is compact and contained in the resolvent, there exists some $\delta_0>0$ such that 
$\{\lambda \in \C: {\rm Re} \lambda \geq -\delta_0, \; {\rm Im}\lambda \leq \lambda_0 \}$ is also in the resolvent. Since we already know that $\Lambda_{\ep,\lambda_0}$ is in the resolvent, it follows that for some $\frac{\pi}{2}>\ep_1 > \ep$ and $0<\delta<\delta_0$ the resolvent set contains $\Sigma_{\ep_1,-\delta}$.  

\qed

Now we denote 
\begin{align}
&AU = \left( \begin{array}{c}
-\vr^*\div \vu\\
\vr_{*}^{-1}[\mu \Delta\vu+(\frac{\mu}{3}+\zeta)\nabla\div\vu+\gamma_{0*}\nabla\eta]  \end{array} \right) \quad {\rm for} \; U=(\eta,\vu) \in \CD_q, \label{bla}\\
& \hat A = A|_{\hat \CD_q}.
\end{align}
Let us consider the Cauchy problem 
\begin{equation} \label{CP1}
\de_t U - \hat A U = 0 \quad {\rm for} \; t>0, \quad  U|_{t=0}=U_0=(\eta_0,\vv_0)\in \hat \CH_q(\Omega) 
\end{equation}
\begin{Proposition} \label{p:decay0}
$\hat A$ generates a $C_0$ semigroup $\{\dot T(t)\}_{t \geq 0}$ on $\hat \CH_q$ and 
\begin{equation} \label{semi2}
\|\dot T(t)U_0\|_{\CH_q} \leq Ce^{-\gamma_1 t}\|U_0\|_{\CH_q}       \end{equation}
for $U_0 \in \hat \CH_q$ and some $\gamma_1>0$. 
\end{Proposition}
{\bf Proof.} The resolvent problem corresponding to \eqref{CP1} is \eqref{res:ME:lin1a}-\eqref{res:bc}. Therefore the generation of analytic semigroup results from the general semigroup theory. More precisely, the semigroup can be expressed by the following formula (see (1.50) in Chapter 9 of \cite{Kato}):   
\begin{equation} \label{form:semi}
\hat T(t) = \frac{1}{2\pi i} \int_{-\infty}^{\infty}e^{\lambda t}(\lambda \I-\hat A)^{-1} {\rm d}\tau,   
\end{equation}
where $\tau = {\rm Im}\,\lambda$ and the line along which we integrate is contained in the resolvent of $\hat A$. Due to Proposition \ref{p5} we can take $\lambda = -\gamma_1+i\tau$ for any $\gamma_1<\delta$ and then the estimate \eqref{semi2} follows 
from \eqref{est:res} and \eqref{form:semi}. 

\qed

Now we can start the proof of Proposition \ref{p2}. Let $\lambda_1>0$ be sufficiently large and $\gamma>0$ sufficiently small to be precised later. For notational convenience let us denote $U=(\eta,\vv)$ let us denote 
\begin{equation} \label{def:Pv}
P_{\vv} U= \vv \quad {\rm for} \, U=(\eta,\vv).
\end{equation}
Now we consider the problem
\begin{equation} \label{prob:U1}
\left\{ \begin{array}{c} 
\de_t U_1 + \lambda_1 U_1 - AU_1=G \quad {\rm in } \; \Omega_0 \times (0,T),\\
P_{\vv} U_1|_{\de \Omega_0}=0, \quad U_1|_{t=0}=U_0.
\end{array} \right.
\end{equation}
Multiplying it by $e^{\gamma t}$ we get
\begin{equation*}
\left\{ \begin{array}{c} 
\de_t (e^{\gamma t}U_1) + (\lambda_1-\gamma)e^{\gamma t}U_1 - Ae^{\gamma t}U_1=e^{\gamma t}G \quad {\rm in }  \; \Omega_0 \times (0,T),\\
P_{\vv}e^{\gamma t}U_1|_{\de \Omega_0}=0, \quad e^{\gamma t}U_1|_{t=0}=U_0.
\end{array} \right.
\end{equation*}
Since we need a maximal regularity estimate for $e^{\gamma t}$ independent from the initial condition, we denote by $G_0$ a zero extension of $G$ on $\R$ and consider a problem
\begin{equation*}
\de_t U_2 + (\lambda_1-\gamma)U_2 - AU_2=e^{\gamma t}G_0 \quad {\rm in} \; \Omega_0 \times \R, \qquad P_{\vv}U_2|_{\de \Omega_0}=0.
\end{equation*}
Applying the Fourier transform in time we obtain 
\begin{equation} \label{prob:U2}
(\lambda_1-\gamma+i\tau) \CF[U_2](\cdot,\tau) - A\CF[U_2](\cdot,\tau)=\CF[e^{\gamma t}G_0](\cdot,\tau) \quad {\rm in} \; \Omega_0, \qquad \CF[P_{\vv}(U_2)]|_{\de \Omega_0}=0.
\end{equation}
Choosing $\lambda_1$ such that $\lambda_1-\gamma>\lambda_0$, where $\lambda_0$ is from Theorem \ref{thm:Rbound} we can deduce 
$\CR$ - bounds analogous to those from Theorem \ref{thm:Rbound}, but for \eqref{prob:U2}. Therefore 
$$
M(\tau):= \CF[U_2](\cdot,\tau) 
$$
is well defined for $\tau \in \mathbb{R}$, $M \in C^1(\R, \CL(\CH_q,\CD_q))$ and moreover \eqref{rbound:Weis} holds. 
Therefore we can apply directly Theorem \ref{thm:Weis} to $M(\cdot)$ to obtain 
\begin{equation} \label{est:U2}
\|\de_t U_2\|_{L_p(\R,\CH_q)}+\|U_2\|_{L_p(\R,D_q)}\leq C\|e^{\gamma t} G_0\|_{L_p(\R,\CH_q)} \leq C\|e^{\gamma t} G\|_{L_p((0,T),\CH_q)}.
\end{equation}
Finally consider the Cauchy problem 
\begin{equation*}
\left\{ \begin{array}{c}
\de_t U_3 + (\lambda_1-\gamma)U_3 - AU_3=0 \quad {\rm in }  \; \Omega \times (0,\infty) \\
P_{\vv}U_3|_{\de \Omega_0}=0, \quad U_3|_{t=0}=U_0-U_2.
\end{array} \right.
\end{equation*}
Let us introduce another constant $\gamma_5>0$.
Then $e^{\gamma_5 t}U_3$ satisfies 
\begin{equation*}
\left\{ \begin{array}{c}
\de_t e^{\gamma_5 t}U_3 + (\lambda_1-\gamma-\gamma_5)e^{\gamma_5 t}U_3 - A(e^{\gamma_5 t}U_3)=0 \quad {\rm in } \;  \Omega \times (0,\infty) \\
P_{\vv}(e^{\gamma_5 t}U_3)|_{\de \Omega_0}=0, \quad e^{\gamma_5 t}U_3|_{t=0}=U_0-U_2.
\end{array}\right.
\end{equation*}
The resolvent problem corresponding to the latter reads 
\begin{equation*}
\left\{ \begin{array}{c}
(\lambda+\lambda_1-\gamma-\gamma_5) e^{\gamma t}U_3 - A(e^{\gamma t}U_3)=0 \quad {\rm in }  \; \Omega \\
P_{\vv}(e^{\gamma t}U_3)|_{\de \Omega_0}=0
\end{array} \right.
\end{equation*}
For $\lambda_1$ sufficiently large we can apply Theorem \ref{thm:Rbound} to obtain the maximal regularity estimate 
\begin{equation} \label{decay0}
e^{\gamma_5 t}\|\de_t U_3\|_{L_p((0,\infty),\CH_q)} 
+ e^{\gamma_5 t}\|U_3\|_{L_p((0,\infty),D_q)} \leq C e^{\gamma_0 t}\|U_0-U_2\|_{W^1_q(\Omega_0)\times B^{2-2/p}_{q,p}(\Omega_0)}
\end{equation}
where $\gamma_0$ is from Theorem \ref{p1}. For sufficiently large $\lambda_1$ we can choose $\gamma_5$ such that $\gamma_5-\gamma_0 \geq \gamma$ and then from \eqref{decay0} we obtain 
\begin{equation} \label{est:U3}
e^{\gamma t}\|\de_t U_3\|_{L_p((0,\infty),\CH_q)} 
+ e^{\gamma t}\|U_3\|_{L_p((0,\infty),D_q)} \leq C \|U_0-U_2\|_{W^1_q(\Omega_0)\times B^{2-2/p}_{q,p}(\Omega_0)}.
\end{equation}
Since $e^{\gamma_t}U_1=U_2+U_3$, from \eqref{est:U2} and \eqref{est:U3} we obtain 
\begin{equation} \label{est:U1}
e^{\gamma t}\|\de_t U_1\|_{L_p((0,\infty),\CH_q)} 
+ e^{\gamma t}\|U_1\|_{L_p((0,\infty),D_q)} \leq C \|\eta_0\|_{W^1_q(\Omega_0)}+\|\vv_0\|_{B^{2-2/p}_{q,p}(\Omega_0)+\|e^{\gamma t} G\|_{L_p((0,T),\CH_q)}}.    
\end{equation}
Now for $U_1=(\eta_1,\vv_1)$ let us define  
\begin{equation} \label{def:tildeU}
\tilde U_1(x,t)=\left(\eta_1(x,t)-\frac{1}{|\Omega_0|}\int_{\Omega_0}\eta_1(y,t)dy, \;\vv_1(x,t)\right).
\end{equation}
Consider the problem 
\begin{equation} \label{prob:V}
\left\{ \begin{array}{c}
\de_t \tilde V - A\tilde V = \lambda_1 \tilde U_1 \quad {\rm in} \; \Omega_0 \times (0,T)\\  
P_{\vv} \tilde V|_{\de \Omega_0}=0, \quad \tilde V|_{t=0}=0.
\end{array} \right.
\end{equation}
Let $\{\hat T(t)\}$ be the semigroup on $\hat \CH_q$ generated by $\hat A$. 
Since $\tilde U_1 \in \hat \CH_q$, we have  
$$
\tilde V = \int_0^t \hat T(t-s)\tilde U_1(\cdot,s)\ds.
$$
Therefore in order to show exponential decay of $\tilde V$ we can use \eqref{semi2} to obtain 
\begin{align*}
&e^{\gamma t}\|\tilde V(\cdot,t)\|_{L_q(\Omega_0)}\leq C \int_0^t e^{\gamma t}e^{-\gamma_1(t-s)}\|\tilde U_1(\cdot,s)\|_{\CH_q}\ds=
C \int_0^t e^{-(\gamma_1-\gamma)(t-s)}e^{\gamma s}\|\tilde U_1(\cdot,s)\|_{\CH_q}\ds\\
&\leq \left(\int_0^t e^{-(\gamma_1-\gamma)(t-s)}\ds\right)^{1/p'}
\left( \int_0^t e^{-(\gamma_1-\gamma)(t-s)}[e^{\gamma s}\|\tilde U_1(\cdot,s)\|_{\CH_q}]^p \right)^{1/p}.
\end{align*}
Therefore 
\begin{align*}
\int_0^T e^{\gamma t}\|\tilde V(\cdot,t)\|_{L_q(\Omega_0)}^p 
\leq C(\gamma_1-\gamma)^{-p/p'}\int_0^T [e^{\gamma s}\|\tilde U_1(\cdot,s)\|_{\CH_q}]^p \ds \int_s^T e^{-(\gamma_1-\gamma)(t-s)} \dt, 
\end{align*}
which implies
\begin{equation} \label{est:tildeV}
\|e^{\gamma t}\tilde V\|_{L_p((0,T),\CH_q)}\leq C(\gamma_1,\gamma,p) \|e^{\gamma t}\tilde U_1\|_{L_p((0,T),\CH_q)}.
\end{equation}
Adding $\lambda_0 \tilde V$ to both sides of \eqref{prob:V} we see that $\tilde V$ satisfies 
\begin{equation}
\left\{ \begin{array}{c}
\de_t \tilde V +\lambda_0 \tilde V - A\tilde V = \lambda_1 \tilde U_1 +\lambda_0 \tilde V \quad {\rm in} \; \Omega_0 \times (0,T)\\  
P_{\vv} \tilde V|_{\de \Omega_0}=0, \quad \tilde V|_{t=0}=0.
\end{array} \right.
\end{equation}
For sufficiently large $\lambda_0$ we have $\CR$ - boundedness for corresponding resolvent problem and from Theorem \ref{thm:maxreg} we get 
\begin{equation*}
\|e^{\gamma t}\de_t \tilde V\|_{L_p(0,T;\CH_q)}+\|e^{\gamma t}\tilde V\|_{L_p(0,T;D_q)} \leq C( \|e^{\gamma t}\tilde U_1\|_{L_p(0,T;D_q)} + \|e^{\gamma t}\tilde V\|_{L_p(0,T;\CH_q)} ),    
\end{equation*}
which combined with \eqref{est:U1} and \eqref{est:tildeV} gives 
\begin{equation} \label{est:V}
\|e^{\gamma t}\de_t \tilde V\|_{L_p(0,T;\CH_q)}+\|e^{\gamma t}\tilde V\|_{L_p(0,T;D_q)} \leq C \left( \|(\eta_0,\vv_0)\|_{W^1_q(\Omega_0)\times B^s_{q,p}(\Omega_0)} + \|e^{\gamma t} G\|_{L_p((0,T),\CH_q)} \right). 
\end{equation}
Next, if we set 
$$
V=\tilde V - \left(\frac{1}{|\Omega_0|}\int_0^t \int_{\Omega_0}\eta(x,t)\dx, 0\right)
$$
then, by \eqref{def:tildeU}, $V$ satisfies 
\begin{equation}
\left\{ \begin{array}{c}
\de_t V - A V = \lambda_1 U_1 \quad {\rm in} \; \Omega_0 \times (0,T)\\  
P_{\vv} \tilde V|_{\de \Omega_0}=0, \quad \tilde V|_{t=0}=0.
\end{array} \right.
\end{equation}
Now if $U_1$ solves \eqref{prob:U1} with $G=(\vf,g)$ then $(\eta,v)=U_1+V$ is the unique solution to \eqref{ME:lin1a} and from \eqref{est:U1} and \eqref{est:V}
we conclude the estimate \eqref{est:decay1} 

\qed

\subsection{Embedding results}
Next we recall some embedding results for Besov spaces.  
The first one is \cite[Theorem 7.34 (c)]{Ad}:
\begin{Lemma} \label{l:imbed1}
Assume $\Omega \in \R^n$ satisfies the cone condition and let $1\leq p,q\leq\infty$ and $sq>n$. Then 
$$
B^s_{q,p}(\Omega_0) \subset C_B(\Omega_0),
$$
where $C_B$ we denote the space of continuous bounded functions.
\end{Lemma}
In particular $u \in B^{2-2/p}_{q,p}(\Omega_0)$ implies $\nabla u \in B^{1-2/p}_{q,p}(\Omega_0)$. Therefore the above Lemma with $s=1-2/p$ yields 
\begin{Corollary} \label{c:imbed1}
Assume $\frac{2}{p}+\frac{3}{q}<1$ and let $\Omega_0$ satisfy the assumptions of Theorem \ref{t1}. Then $B^{2-2/p}_{q,p}(\Omega_0) \subset W^1_{\infty}(\Omega_0)$ and 
\begin{equation}
\|f\|_{W^1_\infty(\Omega_0)} \leq C \|f\|_{B^{2-2/p}_{q,p}(\Omega_0)}.    
\end{equation}
\end{Corollary}
The next result is due to Tanabe  
(cf. \cite[p.10]{Tanabe}):
\begin{Lemma} \label{L:int}
Let $X$ and $Y$ be two Banach spaces such that
$X$ is a dense subset of $Y$ and $X\subset Y$ is continuous.
Then for each $p \in (1, \infty)$  
$$W^1_p((0, \infty), Y) \cap L_p((0, \infty), X) 
\subset C([0, \infty), (X, Y)_{1/p,p})$$
and for every $u\in W^1_p((0, \infty), Y) \cap L_p((0, \infty), X)$ we have
$$\sup_{t \in (0, \infty)}\|u(t)\|_{(X, Y)_{1/p,p}}
\leq (\|u\|_{L_p((0, \infty),X)}^p
+ \|u\|_{W^1_p((0, \infty), Y)}^p)^{1/p}.
$$
\end{Lemma}

\qed 

\section{Local well-posedness}
\subsection{Linearization for the local well-posedness}

Let us start with removing inhomogeneity from the boundary condition \eqref{icbc:lag1}. For this purpose we show 
\begin{Lemma}
Let $\vV$ satisfy the assumptions of Theorem \ref{t1}. Then the problem
\begin{align} \label{def:vub}
&\vr_0\de_t\vu_{b1}-\mu \Delta_y \vu_{b1} - (\frac{\mu}{3}+\zeta)\nabla_y\div_y \vu_{b1} = 0 \qquad \; {\rm in} \; \Omega_0 \times (0,T), \\
&\nonumber \vu_{b1}|_{\Gamma_0}=\vVt, \quad \vu_{b1}|_{t=0}=\vV(0) 
\end{align}
admits a unique solution such that  
\begin{equation} \label{ub:loc}
\|\de_t \vu_{b1}\|_{L_p(0,T;L_q(\Omega_0))}+\|\vu_{b1}\|_{L_p(0,T;W^2_q(\Omega_0))}\leq C \|\de_t\vVt\|_{L_p(0,T;L_q(\Omega_0))}+\|\vVt\|_{L_p(0,T;W^2_q(\Omega_0))}.    
\end{equation}
\end{Lemma}
{\bf Proof.} Denoting $\tilde \vu_{b1} = \vu_{b1}-\vVt$ we have 
\begin{align} \label{def:vuba}
&\vr_0\de_t\tilde \vu_{b1}-\mu \Delta_y \tilde \vu_{b1} - (\frac{\mu}{3}+\zeta)\nabla_y\div_y \tilde \vu_{b1} = \vr_0 \de_t \vVt - \mu \Delta_y \vVt - (\frac{\mu}{3}+\zeta)\nabla_y\div_y \vVt \qquad \; {\rm in} \; \Omega_0 \times (0,T), \\
&\nonumber \tilde \vu_{b1}|_{\Gamma_0}=0, \quad \tilde \vu_{b1}|_{t=0}=0. 
\end{align}
Therefore, if $\vV$ satisfies the assumptions of Theorem \ref{t1} then 
Proposition \ref{p1b} gives 
\begin{equation}
\|\de_t \tilde \vu_{b1}\|_{L_p(0,T;L_q(\Omega_0))}+\|\tilde \vu_{b1}\|_{L_p(0,T;W^2_q(\Omega_0))}\leq C\|\de_t \vVt,\nabla^2_y \vVt\|_{L_p(0,T;L_q(\Omega_0))},     
\end{equation}
which implies \eqref{ub:loc}. 

\qed

As linear system in Theorem \ref{p1} has constant in time coefficients, 
we linearize \eqref{ME:lag1}-\eqref{CE:lag1} around the initial condition. Denoting 
$$
\eta=\vrt-\vr_0, \quad \vv=\vut-\vu_{b1} 
$$
we obtain 
\begin{align} \label{ME:lin1}
&\vr_0 \vv_t - \mu\Delta_y\vv-\left(\frac{\mu}{3}+\zeta \right)\nabla_y\div_y\vv+\gamma_1 \nabla_y \eta= \vF_1(\eta,\vv)\\
\label{CE:lin1}
&\eta_t + \vr_0 \div_y \vv = G_1(\eta,\vv),\\
&\label{icbc:lin1} \vv|_{t=0}=\vu_0-\vV(0), \quad \vv|_{\de \Omega_0}=0.
\end{align}
where $\gamma_1=\pi'(\vr_0)$
and 
\begin{align} \label{F:lin1}
& \vF_1(\eta,\vv) = {\bf R}(\eta+\vr_0,\vv+\vu_{b1}) - \eta \de_t (\vv+\vu_{b1})
-\pi'(\vr_0)\nabla_y \vr_0-[\pi'(\eta+\vr_0)-\pi'(\vr_0)]\nabla_y\eta \\ 
& G_1(\eta,\vv)=G(\eta+\vr_0,\vv+\vu_{b1})-(\eta+\vr_0) \div_y \vu_{b1}-\eta \div_y \vv \label{G:lin1},
\end{align}
and ${\bf R}(\vrt,\vut)$ is defined in \eqref{def:R}.

\subsection{Nonlinear estimates for the local well-posedness.}
We start with a following imbedding result:
\begin{Lemma} \label{l:nonlin1}
Let $f_t \in L_p(0,T;L_q(\Omega_0)),\, f \in L_p(0,T;W^2_q(\Omega_0)) \, f(0,\cdot)\in B^{2-2/p}_{q,p}(\Omega_0)$. Then 
\begin{align} \label{4.1}
{\rm sup}_{t \in (0,T)} \|f\|_{B^{2(1-1/p)}_{q,p}(\Omega_0)} \leq C[ \|f_t\|_{L_p(0,T;L_q(\Omega_0))}+\|f\|_{L_p(0,T;W^2_q(\Omega_0))} +\|f(0)\|_{B^{2-2/p}_{q,p}(\Omega_0)}] \\[5pt] \label{4.2}
{\rm sup}_{t \in (0,T)} \|f\|_{W^1_\infty(\Omega_0))} \leq C[ \|f_t\|_{L_p(0,T;L_q(\Omega_0))}+\|f\|_{L_p(0,T;W^2_q(\Omega_0))} +\|f(0)\|_{B^{2-2/p}_{q,p}(\Omega_0)}].     
\end{align}
\end{Lemma}
{\bf Proof:}
In order to prove \eqref{4.1} we introduce
an extension operator
\begin{equation} \label{def:ext} e_T[f](\cdot, t)
= \begin{cases}
f(\cdot, t) \quad &t\in(0,T), \\
f(\cdot, 2T-t)\quad & t\in(T,2T),\\
0 \quad & t\in(2T,+\infty),
\end{cases}
\end{equation}
If $f|_{t=0}=0$,  then we have
\begin{equation} \label{ext:2} \pd_te_T[f](\cdot, t)
= \begin{cases}
0 \quad &t\in(2T,+\infty),\\
(\pd_tf)(\cdot, t) \quad &t\in(0,T), \\
-(\pd_tf)(\cdot, 2T-t)\quad & t\in(T,2T) \\
\end{cases}
\end{equation}
in a weak sense.
Applying Lemma \ref{L:int} with $X=H^2_q(\Omega_0), \, Y=L_q(\Omega_0)$ and using 
\eqref{def:ext} and \eqref{ext:2} we have
\begin{align*}
&\sup_{t \in (0, T)}\|f(\cdot, t)-f(0)\|_{B^{2(1-1/p)}_{q,p}(\Omega_0)}
\leq \sup_{t \in (0, \infty)}\|e_T[f-f(0)]\|_{B^{2(1-1/p)}_{q,p}(\Omega_0)}\\&\quad 
= (\|e_T[f-f(0)]\|_{L_p((0, \infty), W^2_q(\Omega_0))}^p
+ \|e_T[f-f(0)]\|_{W^1_p((0, \infty), L_q(\Omega_0))}^p)^{1/p}\\
&\quad \leq C(\|f-f(0)\|_{L_p(0, \infty; W^2_q(\Omega_0))}
+ \|\pd_t f\|_{L_p(0, T; L_q(\Omega_0))}).
\end{align*}
This gives \eqref{4.1}, which together with Corollary \ref{c:imbed1} implies \eqref{4.2}.

\qed

Using the results recalled in the previous section and Lemma \ref{l:nonlin1}
we show the following estimate for functions from the space ${\cal Y}(T)$:
\begin{Lemma} \label{l:nonlin2}
Let $(z,\vw) \in B(0,M)\subset {\cal Y}(T)$ and 
$
\|\vw(0)\|_{B^{2-2/p}_{q,p}(\Omega_0)} \leq L.
$
Then 
\begin{align}
&\|\vE^0(\kw),
\nabla_{\kw}\vE^0(\kw),
\|_{L_\infty((0,T)\times\Omega_0)} \leq C(M,L)E(T), \label{est:01}\\
&{\rm sup}_{t \in (0,T)} \|z(\cdot,t)\|_{W^1_q(\Omega_0)}\leq C(M,L)E(T), \label{est:02} \\
&{\rm sup}_{t \in (0,T)}\|\vw(\cdot,t)-\vw(0)\|_{B^{2(1-1/p)}_{q,p}(\Omega_0)}\leq C(M,L),\label{est:03}\\
&\|\vw\|_{L_\infty(0,T;W^1_\infty(\Omega_0))}\leq C(M,L), \label{est:04}
\end{align}
where $\kw$ 
is defined in \eqref{ku}. 
\end{Lemma}
{\bf Proof.} \eqref{est:01} follows immediately from Lemma \ref{l:lag1}.
Next, we have 
$$\|z(\cdot, t)\|_{W^1_q(\Omega_0)}
\leq \int^t_0\|\de_t z(\cdot, s)\|_{W^1_q(\Omega_0)}\,ds
\leq T^{1/{p'}}\|\de_t z\|_{L_p((0, T), W^1_q(\Omega_0))}
\leq C(M)E(T),
$$
which implies \eqref{est:02}. Finally, \eqref{est:03} and \eqref{est:04} follow from \eqref{4.1} and \eqref{4.2}, respectively.  

\qed

Now we can estimate the right hand side of \eqref{ME:lag1} in the regularity required by Theorem \ref{p1}: 
\begin{Lemma} \label{l:nonlin3}
Let $\vF_1(\eta,\vv),G_1(\eta,\vv)$ be defined in \eqref{F:lin1} and \eqref{G:lin1}. Assume that $\vr_0,\vu_0$ and $\vV$ satisfy \eqref{init:norm}. 
Then 
\begin{align} \label{est:FG:loc}
\|\vF_1(\eta,\vv)\|_{L_p(0,T;L_q(\Omega_0))}+\|G_1(\eta,\vv)\|_{L_p(0,T;W^1_q(\Omega_0))}\leq E(T)(\|\eta,\vv\|_{{\cal Y}(T)}+L). 
\end{align}
\end{Lemma}
\bProof 
The proof relies on the estimates collected in Lemma \ref{l:nonlin2}. 
By \eqref{est:02} and \eqref{ub:loc} we have 
\begin{align}
\|\eta \de_t(\vv+\vu_{b1})\|_{\lpq} 
&\leq
\|\eta\|_{L_\infty(\Omega_0 \times (0,T))}\left(\|\de_t\vv\|_{\lpq}+\|\de_t\vu_{b1}\|_{\lpq}\right) \\ \nonumber
&\leq E(T)[\|\eta,\vv\|_{{\cal Y}(T)}+L].    
\end{align}
In order to estimate the remaining terms  
notice that all the quantities \eqref{a2delta}-\eqref{a1div} contain either $\vE(\ku)$ or $\nabla \vE(\ku)$ multiplied by the derivatives of $\vut$ with respect to $y$ of at most second order. Therefore \eqref{est:01} and \eqref{ub:loc} imply 
\begin{equation} \label{est:A} 
\|A_{2\Delta}(\ku)\nabla^2_y\vut,\,
A_{1\Delta}(\ku)\nabla_y\vut,\, 
A_{2\div, i}(\ku)\nabla^2_y\vut,\,
A_{1\div, i}(\ku)\nabla_y\vut\|_{\lpq} \leq E(T)[\|\eta,\vv\|_{{\cal Y}(T)}+L].
\end{equation}
Putting together all above estimates we get the estimate for $\vF_1$. Next, \eqref{est:01} gives immediately 
$$
\|G(\vrt,\vut)\|_{L_p(0,T;W^1_q(\Omega_0))}\leq E(T)[\|\eta,\vv\|_{{\cal Y}(T)}+L], $$
and thus \eqref{est:FG:loc} follows.

\qed

\subsection{Fixed point argument}
Let us define a solution operator 
\begin{equation*}
(\eta,\vv)=S(\bar \eta, \bar \vv) \iff (\eta,\vv) \text{ solves \eqref{ME:lin1}-\eqref{icbc:lin1} with right hand side } \bF_1(\bar \eta, \bar \vv), G_1(\bar \eta, \bar \vv).
\end{equation*}
 By Theorem \ref{p1} and Lemma \ref{l:nonlin3}, $S$ is well defined on ${\cal Y}(T)$ and maps a ball $B(0,M) \subset {\cal Y}(T)$ into itself provided $T$ is sufficiently small w.r.t. $M$ and $L$. Denote 
$$
(\eta_i,\vv_i)=S(\bar \eta_i,\bar \vv_i), \quad i=1,2.  
$$
Then the difference $(\eta_1-\eta_2,\vv_1-\vv_2)$ satisfies 
\begin{align} \label{ME:dif}
&\vr_0 \de_t(\vv_1-\vv_2) - \mu\Delta_y(\vv_1-\vv_2)-\left(\frac{\mu}{3}+\zeta \right)\nabla_y\div_y(\vv_1-\vv_2)+\gamma_1 \nabla_y (\eta_1-\eta_2)= \vF_1(\etab_1,\vvb_1)-\vF_1(\etab_2,\vvb_2)\\
\label{CE:dif}
&\de_t(\eta_1-\eta_2) + \vr_0 \div_y (\vv_1-\vv_2) = G_1(\etab_1,\vvb_1)-G_1(\etab_2,\vvb_2),\\
&\label{icbc:dif} (\vv_1-\vv_2)|_{t=0}=0, \quad (\vv_1-\vv_2)|_{\de \Omega_0}=0,
\end{align}
and we have
\begin{align} \label{dif:F}
&\vF_1(\etab_1,\vvb_1)-\vF_1(\etab_2,\vvb_2)=\bR(\etab_1+\vr_0,\vvb_1+\vu_{b1})-\bR(\etab_2+\vr_0,\vvb_2+\vu_{b1})-(\etab_1-\etab_2)\de_t \vu_{b1}-\etab_1\de_t(\vvb_1-\vvb_2)\\
&\nonumber -\de_t\vvb_2(\etab_1-\etab_2)+\pi'(\vr_0)\nabla_y(\etab_1-\etab_2)-\pi'(\etab_1+\vr_0)\nabla_y(\etab_1-\etab_2)-\nabla_y\etab_2[\pi'(\etab_1+\vr_0)-\pi'(\etab_2+\vr_0)]
\end{align}
and 
\begin{align}\nonumber
&G_1(\etab_1,\vvb_1)-G_1(\etab_2,\vvb_2)=-(\etab_1-\etab_2)\div_y \vu_{b1}-\etab_2\div_y(\vvb_1-\vvb_2)-(\etab_1-\etab_2)\div_y \vvb_1  
-(\etab_1+\vr_0)\vE^1:\nabla_y(\vvb_1-\vvb_2)\\
&-\nabla_y(\vvb_2+\vu_{b1}):[ (\etab_1+\vr_0)(\vE^1-\vE^2)+(\etab_1-\etab_2)\vE^2 ],
\end{align}
where we have denoted 
$$
\vE^1=\vE^0({\bf k}_{\vvb_1+\vu_{b1}}), \quad \vE^2=\vE^0({\bf k}_{\vvb_2+\vu_{b1}}). 
$$
Since $\vE^0(\cdot)$ is smooth, we have 
$$
|\vE^1-\vE^2| \leq C |{\bf k}_{\vvb_1+\vu_{b1}}-{\bf k}_{\vvb_2+\vu_{b1}}|
\leq C {\bf k}_{\vvb_1-\vvb_2}. 
$$
Therefore, recalling the definition of ${\bf R}$ we obtain 
\begin{align*}
\|\bR(\etab_1+\vr_0,\vvb_1+\vu_{b1})-\bR(\etab_2+\vr_0,\vvb_2+\vu_{b1})\|_{L_p(0,T;L_q(\Omega_0))}\leq E(T)\|(\etab_1-\etab_2,\vvb_1-\vvb_2)\|_{{\cal Y}(T)}.
\end{align*}
Estimating the remaining terms on the right hand side of \eqref{dif:F} similarly as in the proof of Lemma \ref{l:nonlin3} we obtain 
\begin{equation} \label{est:dif:F}
\|\vF_1(\etab_1,\vvb_1)-\vF_1(\etab_2,\vvb_2)\|_{L_p(0,T;L_q(\Omega_0))} \leq E(T)\|(\etab_1-\etab_2,\vvb_1-\vvb_2)\|_{{\cal Y}(T)}.
\end{equation}
In a similar way we get 
\begin{equation} \label{est:dif:G}
\|G_1(\etab_1,\vvb_1)-G_2(\etab_2,\vvb_2)\|_{W^1_p(0,T;L_q(\Omega_0))} \leq E(T)\|(\etab_1-\etab_2,\vvb_1-\vvb_2)\|_{{\cal Y}(T)}.    
\end{equation}
Applying \eqref{est:dif:F}, \eqref{est:dif:G} and Theorem \ref{p1} to system \eqref{ME:dif}-\eqref{icbc:dif} we see that $S$ is a contraction on $B(0,M) \subset {\cal Y}(T)$ for sufficiently small times. Therefore it has a unique fixed point $(\eta^*,\vv^*)$. Now 
$$
\vrt= \eta^*+\vr_0, \; \vut=\vv^*+\vu_{b1}
$$
is a solution to \eqref{ME:lag1}-\eqref{icbc:lag1} and 
$$
\|\vrt,\vut\|_{\cal Y(T)}\leq CL.
$$
It is quite standard to verify that after coming back to Eulerian coordinates we obtain a solution with the estimate \eqref{est:t1}, however for the sake of completeness we justify it briefly in the next subsection. 

\subsection{Equivalence of norms in Lagrangian and Eulerian coordinates} \label{sec:equiv}
By \eqref{est:04}, the Jacobian of the transformation $\bXu$ is bounded in space-time. Therefore, Lemma \ref{l:lag2} implies the equivalence of $L_p(0,T;L_q)$ norms of a function and its first-order space derivatives. Furthermore, we have  
$$
\nabla^2_y \tilde f(t,y)=\nabla_xf(t,\bXu(t,y))\nabla^2_y \bXu + (\nabla_y \bXu)^2 \nabla^2_x f(t,\bXu(t,y)).  
$$
Again by \eqref{est:04}, $\nabla_y \bXu$ is bounded in space-time, which together with embedding $W^1_q(\Omega_t) \subset L_\infty(\Omega_t)$ for $t \in [0,T)$ gives equivalence of $L_p(L_q)$ norms of second space derivatives. However, we have a different situation for the time derivative. 
The solution constructed in Lagrangian coordinates satisfies 
$$
\vrt_t \in L_p(0,T;W^1_q(\Omega_0)). 
$$
However, due to \eqref{dt:lag} this does not imply the same regularity for the density in Eulerian coordinates. Nevertheless, the regularity of $\vu$ implies  
$$
\vr_t \in L_p(0,T;L_q(\Omega_t)), 
$$
which is the regularity in the assertion of Theorem \ref{t1}.

\section{Global well-posedness}

\subsection{Linearization}
Again we first reduce the problem to homogeneous boundary condition. 
\begin{Lemma}
If $\vV$ satisfies the assumptions of Theorem \ref{t2} then 
the problem
\begin{align} \label{def:vub2}
&\vr^*\de_t\vu_{b2}-\mu \Delta_y \vu_{b2} - (\frac{\mu}{3}+\zeta)\nabla_y\div_y \vu_{b2} = 0  \qquad \; {\rm in} \; \Omega_0 \times (0,T), \\
&\nonumber \vu_{b2}|_{\Gamma_0}=\vVt, \quad \vu_{b2}|_{t=0}=\vV(0) 
\end{align}
admits a unique global in time solutions
$\vu_{b2}$ with the decay estimate 
\begin{equation} \label{ub:glob}
\|e^{\gamma t}\de_t \vu_{b2}\|_{L_p(0,T;L_q(\Omega_0))}+\|e^{\gamma t}\nabla_y \vu_{b2}\|_{L_p(0,T;W^1_q(\Omega_0))}+\|e^{\gamma t}(\vu_{b2}-\vVt))\|_{L_p(0,T;L_q(\Omega_0))}\leq C \|e^{\gamma t}(\de_t \vVt,\nabla^2_y \vVt)\|_{L_p(0,T;L_q(\Omega_0))}.    
\end{equation}
\end{Lemma}
{\bf Proof.} Let us define $\tilde \vu_{b_2}=\vu_{b2}-\vVt$.
If $\vV$ satisfies the assumptions of Theorem \ref{t2} then Theorem \ref{p2} implies
\begin{equation*}
\|e^{\gamma t}\de_t \tilde \vu_{b2}\|_{L_p(0,T;L_q(\Omega_0))}+\|e^{\gamma t} \tilde \vu_{b2}\|_{L_p(0,T;W^2_q(\Omega_0))}\leq C\|e^{\gamma t}(\de_t \vVt,\nabla^2_y \vVt)\|_{L_p(0,T;L_q(\Omega_0))},     
\end{equation*}
which gives \eqref{ub:glob}.

\qed

This time we have to linearize the density around the constant $\vr^*$. Denoting 
$$\sigma=\vrt-\vr^*, \quad \vv=\vut-\vu_{b2} $$ 
we obtain from \eqref{ME:lag1}-\eqref{icbc:lag1} 
\begin{align} \label{ME:lin2}
&\vr^* \vv_t - \mu\Delta_y\vv-\left(\frac{\mu}{3}+\zeta \right)\nabla_y\div_y\vv+\gamma_2 \nabla_y \sigma= \vF_2(\sigma,\vv),\\
\label{CE:lin2}
&\sigma_t + \vr^* \div_y \vv = G_2(\sigma,\vv),\\
&\label{icbc:lin2} \vv|_{t=0}=\vu_0-\vV(0), \quad \vv|_{\de \Omega_0}=0.
\end{align}
where $\gamma_2=\pi'(\vr^*)$
and 
\begin{align} \label{F:lin2}
& \vF_2(\sigma,\vv) = \vF(\sigma+\vr^*,\vv+\vu_{b2}) - \sigma \de_t (\vv+\vu_{b2})
-[\pi'(\sigma+\vr^*)-\pi'(\vr^*)]\nabla_y\sigma, \\ 
& G_2(\sigma,\vv)=G(\sigma+\vr^*,\vv+\vu_{b2})-(\sigma+\vr^*) \div_y \vu_{b2}-\sigma \div_y \vv . \label{G:lin2}
\end{align}

\subsection{Nonlinear estimates for the global well-posedness}
We start with an analog of Lemma \ref{l:nonlin2} which will be used to estimate the nonlinearities for large times. We also recall the definitions \eqref{def:Y} and \eqref{def:dotY} of norm and seminorm on ${\cal Y}(\infty)$.
\begin{Lemma} \label{l:nonlin4}
Let $e^{\gamma t}(z,\vw)\in {\cal Y}(\infty)$ for some $\gamma>0$ and $$
\|z(0)-\vr^*\|_{W^1_p(\Omega_0)}+\|\vw(0)\|_{B^{2-2/p}_{q,p}(\Omega_0)}\leq \ep.
$$
Then 
\begin{align}
&\|\vE^0(\kw),
\nabla_{\kw}\vE^0(\kw),
\|_{L_\infty((0,\infty)\times\Omega_0)} \leq C\expnormzw, \label{est:05}\\
&{\rm sup}_{t \in (0,\infty)} \|z(\cdot,t)\|_{W^1_q(\Omega_0)}\leq C[\epsilon+\expnormzw] , \label{est:06} \\
&{\rm sup}_{t \in (0,\infty)}\|\vw(\cdot,t)-\vw(0)\|_{B^{2(1-1/p)}_{q,p}}\leq C\expnormzw, \label{est:07}\\
&\|\vw\|_{L_\infty(0,\infty,W^1_\infty(\Omega_0))}\leq C\expnormzw, \label{est:08}
\end{align}
where $\kv$ 
is defined in \eqref{ku}. 
\end{Lemma}
{\bf Proof:} We have 
$$
\int_0^\infty \|\nabla_y \vw\|_{\infty} \dt\leq \left(\int_0^\infty e^{-\gamma t p'} \dt\right)^{1/p'} \left(\int_0^\infty e^{\gamma t p}\|\vw\|_{W^2_q(\Omega_0)} \dt\right)^{1/p},
$$
which implies \eqref{est:05}. Next, 
$$
\|\sigma(\cdot,t)\|_{L_\infty(\Omega_0)}\leq \|z(0)-\vr^*\|_{L_\infty(\Omega_0)}+\int_0^t\|z_t(s,\cdot)\|_{L_\infty(\Omega_0)} \dt
$$$$
\leq \epsilon+ C\left(\int_0^t e^{-\gamma s p'} \d s\right)^{1/p'} \left(\int_0^\infty e^{\gamma s p}\|z_t\|_{W^1_q(\Omega_0)} \d s\right)^{1/p},
$$
which yields \eqref{est:06}. 
Finally, \eqref{est:07} and \eqref{est:08} follows from Lemma \ref{l:nonlin1}. 

\qed

The following lemma gives estimates for the right hand sides of \eqref{ME:lin2}-\eqref{CE:lin2}.
\begin{Lemma} \label{l:nonlin5}
Let $\vF_2(\vrt,\vut),G_2(\vrt,\vut)$ be defined in \eqref{F:lin2} and \eqref{G:lin2}. Assume that $\vr_0,\vu_0$ and $\vV$ satisfy \eqref{init:norm2}. 
Then 
\begin{align} \label{est:F}
\|\vF_2(\vrt,\vut)\|_{L_p(0,\infty;L_q(\Omega_0))}+\|G_2(\vrt,\vut)\|_{L_p(0,\infty;W^1_q(\Omega_0))}\leq C(\expnorm^2+\epsilon). 
\end{align}
\end{Lemma}
{\bf Proof}. First, analogously to \eqref{est:A}, this time using \eqref{est:05}
and \eqref{ub:glob} we obtain 
\begin{align} \label{est:A2}
&\|e^{\gamma t}A_{2\Delta}(\ku)\nabla^2_y\vut,\,
A_{1\Delta}(\ku)\nabla_y\vut,\, 
A_{2\div, i}(\ku)\nabla^2_y\vut,\,
A_{1\div, i}(\ku)\nabla_y\vut\|_{L_p(0,\infty;L_q(\Omega_0))} \nonumber\\ 
&\leq C \expnorm \|e^{\gamma t}(\vv+\vu_{b2})\|_{L_p(0,\infty;L_q(\Omega_0))}\nonumber\\ 
&\leq C \expnorm[\expnorm + \|e^{\gamma t}(\de_t\vVt,\nabla^2_y\vVt)\|_{L_p(0,\infty;L_q(\Omega_0))}] \nonumber \\
&\leq C \expnorm[\expnorm +\epsilon].
\end{align}
Next, by \eqref{est:06} and \eqref{ub:glob} 
$$
\|e^{\gamma t}\de_t(\vv+\vu_{b2})\|_{L_p(0,\infty;L_q(\Omega_0))}\leq \|\sigma\|_{L_\infty((0,\infty)\times\Omega_0)}\|e^{\gamma t}\de_t(\vv+\vu_{b2})\|_{L_p(0,\infty;L_q(\Omega_0))} 
\leq C[\ep+\expnorm]^2,
$$
and 
$$
\|e^{\gamma t}[\pi'(\vrt)-\pi'(\vr^*)]\nabla_y \sigma\|_{L_p(0,\infty;L_q(\Omega_0))} 
\leq \|\sigma\|_{L_\infty((0,\infty)\times \Omega_0)}\|e^{\gamma t}\nabla_y \sigma\|_{L_p(0,\infty;L_q(\Omega_0))} \leq C[\ep+\expnorm]^2.
$$
Combining all above estimates we get the required estimate for $\|\vF_2\|_{L_p(0,\infty;L_q(\Omega_0))}$. Finally, $G_2$ and its space derivatives 
are estimated in a similar way using Lemma \ref{l:nonlin4} and \eqref{ub:glob}.

\qed

\subsection{Proof of Theorem \ref{t2}}
It is now easy to verify the following estimate which allows to prolong the local solution for arbitrarily large times.
\begin{Lemma} Assume $\sigma,\vv$ is solution to \eqref{ME:lin2}-\eqref{icbc:lin2} with $\vr_0,\vu_0$ and $\vV$ satisfying the assumptions of Theorem \ref{t2}. Then 
\begin{equation} \label{est:expnorm}
\expnorm \leq E(\epsilon).    
\end{equation}
\end{Lemma}
{\bf Proof.} Combining Theorem \ref{p2} and Lemma \ref{l:nonlin5} we obtain 
\begin{equation} \label{est:expnorm1}
\|e^{\gamma t} (\sigma,\vv)\|_{\dot {\mathcal{Y}}(T)} \leq C[\epsilon + \|e^{\gamma t} (\sigma,\vv)\|_{\dot {\mathcal{Y}}(T)}^2].   
\end{equation}
Note that we derived this inequality for $T = \infty$, however it is easy to observe that the same arguments yield \eqref{est:expnorm1} for any $T > 0$. Consider the equation 
$$
x^2-\frac{x}{C}+\epsilon =0. 
$$
Its roots are 
$$
x_1(\epsilon)= \frac{1}{2C}-\sqrt{\frac{1}{4C^2}-\epsilon},
\quad x_2(\epsilon)= \frac{1}{2C}+\sqrt{\frac{1}{4C^2}-\epsilon}. 
$$
Notice that the inequality \eqref{est:expnorm1} implies either $\|e^{\gamma t} (\sigma,\vv)\|_{\dot {\mathcal{Y}}(T)} \leq x_1(\epsilon)$
or $\|e^{\gamma t} (\sigma,\vv)\|_{\dot {\mathcal{Y}}(T)} \geq x_2(\epsilon)$. However, 
$$
\|e^{\gamma t} (\sigma,\vv)\|_{\dot {\mathcal{Y}}(T)} \rightarrow 0
$$
as $T \to 0$, therefore 
$$
\|e^{\gamma t} (\sigma,\vv)\|_{\dot {\mathcal{Y}}(T)} \leq x_1(\epsilon)
$$
for $T$ small. Finally, $\|e^{\gamma t} (\sigma,\vv)\|_{\dot {\mathcal{Y}}(T)}$ is continuous in time and therefore
$$
\|e^{\gamma t} (\sigma,\vv)\|_{\dot {\mathcal{Y}}(\infty)} \leq x_1(\epsilon).
$$

\qed

Now it is a standard matter to prolong the local solution for arbitrarily large times. For this purpose it is enough to observe that if the initial data satisfies the smallness assumption from Theorem \ref{t2} then the time of existence from Theorem \ref{t1} satisfies $T>C(\epsilon)>0$. Therefore, 
for arbitrarily large $T^*$ we can obtain a solution on $(0,T^*)$ in a finite number of steps. By the estimate \eqref{est:expnorm} this solution satisfies \eqref{est:glob1}-\eqref{est:glob2}.  

Finally, the equivalence of norms can be justified as in Section \ref{sec:equiv}, using \eqref{est:08} instead of \eqref{est:04}.

\section*{Acknowledgements}
 The works of O. Kreml, \v S. Ne\v casov\'a,  T.Piasecki were supported by the Czech Science Foundation grant GA19-04243S in the framework of RVO 67985840. The work of T. Piasecki was partially supported by Polish National Science Centre grant 2018/29/B/ST1/00339.

\end{document}